\documentclass[11pt]{amsart}

\usepackage{amssymb,amscd,amsthm,amsxtra}
\usepackage{latexsym}
\usepackage[usenames]{color}

\newtheorem{thm}{Theorem}[section]
\newtheorem{cor}[thm]{Corollary}
\newtheorem{lem}[thm]{Lemma}
\newtheorem{prop}[thm]{Proposition}
\theoremstyle{definition}
\newcommand{\comment}[1]{}

\theoremstyle{remark}
\newtheorem{rem}[thm]{Remark}
\numberwithin{equation}{section}
\newcommand{\R}{\mathbb R}
\newcommand{\C}{\mathbb C}

\newcommand{\CN}{\mathcal N}

\begin{document}

\title[Global well-posedness for defocusing $L^{2}$-critical NLS in 1d]
{\bf Global well-posedness and polynomial bounds for the defocusing
$L^{2}$-critical nonlinear Schr\"odinger equation in $\R$}

\author{Daniela De Silva}
\address{Department of Mathematics, Johns Hopkins University, Baltimore, MD 21218}
\email{\tt  desilva@math.jhu.edu}

\author{Nata\v{s}a Pavlovi\'{c}}
\address{Department of Mathematics, Princeton University, Princeton, NJ 08544-1000}
\email{\tt natasa@math.princeton.edu}
\thanks{N.P. was supported by N.S.F. Grant DMS 0304594.}

\author{Gigliola Staffilani}
\address{Department of Mathematics, Massachusetts Institute of Technology, Cambridge, MA 02139-4307} \email{\tt gigliola@math.mit.edu}
\thanks{G.S. was supported by N.S.F. Grant DMS 0602678.}

\author{Nikolaos Tzirakis}
\address{Department of Mathematics, University of Toronto, Toronto, Ontario, Canada M5S 2E4}
\email{tzirakis@math.toronto.edu}

\date{March 11, 2007}


\begin{abstract}
 We prove global well-posedness for low regularity data for the
one dimensional quintic defocusing nonlinear Schr\"odinger
equation. Precisely we show that a unique and global solution
exists for initial data in the Sobolev space $H^{s}(\mathbb R)$
for any $s>\frac{1}{3}$. This improves the result in \cite{tz},
where global well-posedness was established for any
$s>\frac{4}{9}$. We use the $I$-method to take advantage of the
conservation laws of the equation. The new ingredient in our proof
is an interaction Morawetz estimate for the smoothed out solution
$Iu$. As a byproduct of our proof we also obtain that the $H^{s}$
norm of the solution obeys polynomial-in-time bounds.

\end{abstract}
\maketitle

\section{Introduction}

In this paper we study the global well-posedness of the following
 initial value problem (IVP) for the $L^2$-critical defocusing nonlinear
Schr\"{o}dinger equation (NLS): \vspace{1mm}
\begin{align}
&iu_{t}+ \Delta u -|u|^{4}u=0, \label{ivp1}\\
&u(x,0)=u_{0}(x)\in H^{s}({\mathbb R}), \; \; \; \; x \in {\mathbb R},\;  t\in {\mathbb
R},
\label{bc1}
\end{align}
where  $H^s$ denotes the usual inhomogeneous Sobolev space of order $s$.

 We adopt the standard notion of local well-posedness,
that is, we say that the IVP above  is locally well-posed in
$H^s$ if for any initial data $u_0 \in H^s$ there exists a
positive time $T = T(\|u_0\|_{H^{s}})$ depending only on the norm
of the initial data, such that a solution to the initial value
problem exists on the time interval $[0,T]$, is unique in a
certain Banach space of functions $X\subset C([0,T],H^s_{x})$, and
the solution map from $H^s_{x}$ to $C([0,T],H^s_{x})$ depends
continuously on the initial data on the time interval $[0,T]$. If
$T$ can be taken arbitrarily large we say that the IVP is globally well-posed.

Local well-posedness for the initial value problem
\eqref{ivp1}-\eqref{bc1} in $H^s$ for any $s > 0$ was established
in \cite{cw}, see also \cite{cw1}. A local solution  also exists
for $L^{2}$ initial data \cite{cw1}, but the time of existence
depends not only on the $H^{s}$ norm of the initial data, but also
on the profile of $u_{0}$. For more details on local existence
see, for example, \cite{jb2}, \cite{cw1} and \cite{tt}.

Local in time solutions of \eqref{ivp1}-\eqref{bc1} enjoy mass
conservation
\begin{equation}\label{mass}
\|u(\cdot,t)\|_{L^2(\mathbb R^n)} =
\|u_0(\cdot)\|_{L^2(\mathbb R^n)}.
\end{equation}
Moreover, $H^{1}$ solutions enjoy conservation of the energy
\begin{equation}
E(u)(t)=\frac{1}{2}\int_{\mathbb R} |\partial_{x} u(t)|^{2}dx+\frac{1}{6}\int_{\mathbb R}
|u(t)|^{6}dx=E(u)(0),
\end{equation}
which together with \eqref{mass} and the local theory immediately
yields global in time well-posedness for \eqref{ivp1}-\eqref{bc1}
with initial data in $H^1$.  On the other hand the $L^{2}$ local well-posedness and
the conservation of the $L^{2}$ norm do not immediately imply
global well-posedness as in the case of the finite energy data since in
this case $T=T(u_{0})$. It is conjectured though that
\eqref{ivp1}-\eqref{bc1} is globally well-posed for initial data
in $H^s$, $s \geq 0$.

Existence of global solutions to \eqref{ivp1}-\eqref{bc1}
corresponding to initial data below the energy treshold was first
obtained in \cite{ckstt1}, \cite{ckstt5} by using the method of
``almost conservation laws", or ``$I$-method" (for a detailed
description of this method see section \ref{sec-I} below). The gap
between known local and global well-posedness was further filled
out  in \cite{tz},   where  global well-posedness was obtained in
$H^s({\mathbb R})$ with $s > 4/9$.  The question of global
well-posedness below the energy space $H^{1}$  we described above
can also be formulated in $\R^{d}$, where the equivalent $L^{2}$
critical problem has nonlinearity $|u|^{4/d}u$. For $d=2$ the
$I$-method was  also partially successful in \cite{ckstt2}, and in
\cite{fg}, where global well-posedness was obtained in $H^{s}$ for
$s\geq 1/2$ . This last paper is particularly interesting since it
combines the $I$-method with a local in time Morawetz type
estimate. In \cite{dpst-high} this idea was then extended to all
dimensions $d\geq 3$. Still none of these results reached the
space $L^{2}$. Recently  in \cite{tvz} the authors proved global
well-posedness in $L^2$ and scattering for the $L^2$-critical NLS
problem in all dimensions $d \geq 3$, assuming spherically
symmetric initial data.  The proof relies upon a combination of
several sophisticated tools among which compensated compactness
and a frequency-localized Morawetz estimate. These tools seem not
to be enough  when $d=1,2$ or when the  data are no longer radial.

In this paper we only consider the case $d=1$ and we  prove the following result:

\begin{thm}\label{th-gwp}
The initial value problem \eqref{ivp1}-\eqref{bc1} is globally
well-posed in $H^{s}(\mathbb R)$, for any $1>s>\frac{1}{3}$. Moreover the
solution satisfies
$$\sup_{t \in [0,T]}\|u(t)\|_{H^{s}(\mathbb R)} \leq C(1+T)^{\frac{s(1-s)}{2(3s-1)}}$$
where the constant $C$ depends only on $s$ and
$\|u_{0}\|_{L^{2}}$.
\end{thm}

We prove Theorem \ref{th-gwp} by combining the $I$-method with an
interaction Morawetz-type estimate for the smoothed out version
$Iu$ of the solution. Such a Morawetz estimate for an almost
solution, that below we call  ``almost Morawetz'',
is the main novelty of this paper. As mentioned above, the approach of
combining the $I$-method with an interaction Morawetz estimate was
described in \cite{fg} where the
$L^2$-critical NLS in 2d  was treated; see also \cite{ckstt4}.
However in order to
obtain global well-posedness for \eqref{ivp1}-\eqref{bc1}
corresponding to initial data in $H^s$ with $s \leq 4/9$\footnote{
Note that global well-posedness for initial data in $H^s$ with $s
> 4/9$ was established in \cite{tz}.} it was not
enough to obtain an interaction Morawetz estimate for the solution
$u$ itself (as in the case of \cite{ckstt4}, \cite{dpst-high},
\cite{fg}).

Before giving an outline of the proof we say a few words about the
above mentioned tools: the $I$-method and the ``almost Morawetz''
estimate.

The $I$-method was first  introduced  by  Colliander et al  (see, for example,
\cite{ckstt1, ckstt2, ckstt4}).
It is based on the almost conservation of a certain modified
energy functional. The idea is to replace the conserved quantity
$E(u)$, which is no longer available for $s<1$, with an ``almost
conserved'' variant $E(Iu)$, where $I$ is a smoothing operator of
order $1-s$, which behaves like the identity for low frequencies
and like a fractional integral operator for high frequencies.
However $Iu$ is not a solution to \eqref{ivp1} and hence one
expects an energy increment. This increment is quantifying $E(Iu)$
as an ``almost conserved'' energy. The key is to prove that on
intervals of fixed length, where local well-posedness is
satisfied, the increment of the modified energy $E(Iu)$ decays
with respect to a large parameter $N$ (for the precise definition
of $I$ and $N$ we refer the reader to Section \ref{sec-I}). This
requires delicate estimates on the commutator between $I$ and the
nonlinearity. When $d=1$,  hence the nonlinearity is algebraic,
one can write the commutator explicitly using the Fourier
transform, and control it by multi-linear analysis and bilinear
estimates. The analysis above can be carried out in the $X^{s,b}$
spaces setting, where one can use the smoothing bilinear
Strichartz  estimate of Bourgain (see e.g. \cite{jb1}) along with
Strichartz estimates, to demonstrate the existence of global rough
solutions (see \cite{ckstt1, ckstt5} and \cite{tz}).

We now turn to our second tool: the ``almost Morawetz estimate", that is an \textit{a priori}
interaction Morawetz-type estimate for the ``approximate
solution'' $Iu$ to the initial value problem
\begin{align}
&iIu_{t}+ \Delta Iu -I(|u|^{4}u)=0, \label{aivp1}\\
&Iu(x,0)=Iu_{0}(x)\in H^{1}({\mathbb R}),\; \; \; \; x \in {\mathbb R} , t\in {\mathbb R}
\label{abc1}.
\end{align}
For the original problem \eqref{ivp1} one can prove that solutions
satisfy the following a priori bound (see \cite{chvz06})
$$\|u\|_{L_{t \in [0,T]}^{8}L_{x}^{8}}^{8}
\lesssim \sup_{t \in [0,T]}\|u\|_{\dot{H}^{1}}\|u\|_{L^{2}}^{7}.$$
For initial data below $H^{1}$ this estimate is not useful anymore
since $u$ is not in $H^{1}$. We introduce the $I$-operator with
the aim of getting an a priori estimate of the form
$$\|Iu\|_{L_{t \in [0,T]}^{8}L_{x}^{8}}^{8}
\lesssim \sup_{t \in [0,T]}\|Iu\|_{\dot{H}^{1}}\|Iu\|_{L^{2}}^{7}+
Error,$$ where the $Error$ terms are negligible in some sense. To
achieve this, we work with mixed Lebesgue spaces that we denote by
$S_{I}(J)$ and define\footnote{See Section 2 for a definition of
the operator $\langle\partial_{x}\rangle $.} as
$$\label{SI} S_I(J) := \{f\; |\; \sup_{(q,r) \ \ admissible}
\|\langle\partial_{x}\rangle If\|_{L^q_tL^r_x(J\times \mathbb R)} < \infty\},$$
where the pair $(q,r)$ is said to be admissible if $\frac{2}{q}+\frac{1}{r}=\frac{1}{2}$
and $2\leq q,r\leq \infty$.

 We show that inside the interval where the local solution
exists, the error term is very small. The proof of this fact
relies on harmonic analysis estimates of Coifman-Meyer type and is
given in Section \ref{sec-almor}\footnote{Here we were not able
to use the $X^{s,b}$ spaces machinery. However the approach that
we pursue creates no additional difficulties, since we can iterate
our solutions in the intersection of the spaces $S_I$ and
$X^{s,b}$, see Proposition \ref{lwp} for details.}. Because of the fact that we work with local-in-time solutions
 we restrict the above a priori bounds to local intervals of the form $[t_{0},t_{1}]$.

\vspace{2mm}

 Now we give an outline of the proof to show how to combine our two tools. Fix a large value of
time $T_{0}$. If $u$ is a solution to \eqref{ivp1} in the time
interval $[0,T_{0}],$ then
$u^\lambda(x)=\frac{1}{\lambda^{\frac{1}{2}}}u(\frac{x}{\lambda},\frac{t}{\lambda^{2}})$
is a solution to the same equation in $[0, \lambda^{2}T_{0}]$. We
choose the parameter $\lambda>0$ so that
$E(Iu_{0}^{\lambda})=O(1)$. Using Strichartz estimates we show
that if $J=[t_{0},t_{1}]$ and $\|Iu^{\lambda}\|_{L^{6}_tL^{6}_x(J\times
\Bbb R^n )}^{6}<\mu$, where $\mu$ is a small universal constant,
then for $b$ close to $1/2$
$$ \|u^{\lambda}\|_{S_I(J)} \lesssim \|Iu^\lambda(t_{0})\|_{H^1},$$
and
$$\|u^{\lambda}\|_{X^{s,b}_I(J)} \lesssim \|Iu^\lambda(t_{0})\|_{H^1},$$
where
$$X_{I}^{s,b}(J) = \{ f | \; \; \|If\|_{X^{s,b}(J \times {\mathbb R})} < \infty \}.$$
Moreover in  this same time interval where the problem is
well-posed, we can prove the following ``almost conservation law'', provided
that\footnote{We use the second modified
energy approach as in \cite{ckstt5} and \cite{tz}.}
$s >1/4$.
\begin{equation}\label{acl}
|E^{2}(u^{\lambda})(t_{1}) - E^{2}(u^{\lambda})(t_{0})|
\lesssim N^{-2}\|Iu^{\lambda}\|_{X^{s,b}([t_{0},t_{1}]\times \mathbb R)}^{10}
\lesssim N^{-2}\|Iu^\lambda(t_{0})\|_{H^1}^{10} \lesssim N^{-2},
\end{equation} 
where $E^{2}(u^{\lambda})$ denotes the second modified energy (for the definition 
see Section 3). 
Of course, for the arbitrarily large interval $[0,\lambda^{2}T_{0}] $
we do not have
$$\|Iu^{\lambda}\|_{L^{6}_tL^{6}_x([0,\lambda^{2}T_{0}]\times
\mathbb R^n )}^{6}<\mu.$$ However we can interpolate the a priori
information that we have about the $L^8_tL^8_x$ norm of $Iu$ with
the following a priori bound
$$\|Iu^{\lambda}\|_{L^{2}_{t\in [0,T]}L^{2}_x}
\leq T^{\frac{1}{2}} \|Iu^{\lambda}\|_{L^{\infty}_{t\in
[0,T]}L^{2}_x } \leq T^{\frac{1}{2}}
\|u^{\lambda}\|_{L^{\infty}_{t \in [0,T]}L^{2}_x}
=T^{\frac{1}{2}}\|u_{0}^{\lambda}\|_{L^{2}_x}=
T^{\frac{1}{2}}\|u_{0}\|_{L^{2}_x},$$ obtained using H\"older's
inequality, the definition of the $I$, mass conservation and the
fact that the problem is  $L^{2}$-critical. Then we get an
$L_{t}^6L_{x}^6$ bound valid for $J=[t_{0},t_{1}]$ and we use this bound
to partition the arbitrarily large interval $[0,\lambda^{2}T_{0}]$
into $L$ intervals where the local theory uniformly applies.
$L=L(N,T)$ is finite and defines the number of the intervals in
the partition that will make the Strichartz $L^{6}_tL^{6}_x$ norm
of $Iu$ less than $\mu$ in each interval. The next step is to prove that the second
modified energy $E^{2}(u^{\lambda})$ is just an approximation of
the first modified energy $E^1(u^{\lambda}) = E(Iu^{\lambda})$ in the sense that
\begin{equation} \label{intropert}
E^{2}(u^{\lambda}) \sim E^{1}(u^{\lambda}) + O(\frac{1}{N^{\epsilon}}),
\end{equation}
for $s > 1/3$ and $N \gg 1$. Since $E(Iu^{\lambda})$ controls the $H^1$ norm of $Iu$ we have by (\ref{acl})
$$\|Iu^{\lambda}\|_{H^{1}} \lesssim LN^{-2}.$$
To maintain the bound $\|Iu^{\lambda}\|_{H^1} \lesssim 1$ we must
have that
$$L(N,T)\sim N^{2},$$ and this condition will require $s>1/3$. 
One should notice already that  the restriction  $s>1/3$ 
appears in two separate parts of our proof: in \eqref{intropert} and  
in the last step recorded above. So this regularity is a threshold to the method, at least if
 one wants to use the second modified energy.
For a more detailed proof the reader should check Section \ref{sec-proof}.

We conclude this introduction by announcing that work in progress of the fourth author in collaboration with
J. Colliander and M. Grillakis, shows that a similar approach
can be used in the $L^{2}$-critical case when $d=2$.

\subsection*{Organization of the paper}
In Section 2 we introduce some notation and state important
propositions that we will use throughout the paper. In Section 3
we review the $I$ method, prove the local well-posedness theory
for $Iu$ and obtain an upper bound on the increment of the second
modified energy. In Section 4 we prove the ``almost Morawetz''
inequality which is the heart of our argument. Finally in Section
5 we give the details of the proof of global well-posedness stated
in Theorem \ref{th-gwp}.
\subsection*{Acknowledgements} We would like to thank Jim Colliander for helpful conversations and suggestions.

\section{Notation and preliminaries}

\subsection {Notation}
In what follows we use $A \lesssim B$ to denote an estimate of the
form $A\leq CB$ for some constant $C$. If $A \lesssim B$ and $B
\lesssim A$ we say that $A \sim B$. We write $A \ll B$ to denote
an estimate of the form $A \leq cB$ for some small constant $c>0$.
In addition $\langle a \rangle:=1+|a|$ and $a\pm:=a\pm \epsilon$
with $0 < \epsilon <<1$.  The
reader also has to be alert that we sometimes do not explicitly
write down constants that depend on the $L^{2}$ norm of the
solution. This is justified by the conservation of the $L^{2}$
norm.

\subsection {Definition of spaces}
We use $L^r_x(\mathbb R)$ to denote the Lebesgue space of
functions $f : \R \rightarrow \C$ whose norm
$$\|f\|_{L^r_x}:=\left( \int_{\R^n}|f(x)|^r dx
\right)^{\frac{1}{r}}$$ is finite, with the usual modification in
the case $r=\infty.$ We also use the mixed space-time Lebesgue
spaces $L^q_tL^r_x$ which are equipped with the norm
$$\|u\|_{L^q_tL^r_x} := \left(\int_J \|u\|_{L^r_x}^q dt\right)^{\frac{1}{q}}$$
for any space-time slab $J \times \mathbb R,$ with the usual modification when
either $q$ or $r$ are infinity. When $q=r$ we abbreviate
$L^q_tL^r_x$ by $L^q_{t,x}.$

As usual, we define the Fourier transform of $f(x) \in L_{x}^{1}$ by
$$\hat{f}(\xi)=\int_{\mathbb R} e^{-2\pi i\xi x}f(x)dx.$$
For an appropriate class of functions the following Fourier inversion formula holds:
$$f(x)=\int_{\mathbb R} e^{2\pi i\xi x}\hat{f}(\xi)d\xi.$$
We define the fractional differentiation operator
$|\partial_{x}|^\alpha$ for any real $\alpha$ by
$$\widehat{|\partial_{x}|^\alpha u}(\xi):=
|\xi|^\alpha \hat{u}(\xi),$$
and analogously
\begin{equation}\label{operator}
    \widehat{\langle\partial_{x}\rangle^\alpha u}(\xi):=
\langle\xi \rangle^\alpha \hat{u}(\xi).$$
The inhomogeneous Sobolev space $H^s$ is given via
$$\|u\|_{H^s} = \|\langle\partial_{x}\rangle^s
u\|_{L^2_x},
\end{equation}
while the homogeneous Sobolev space $\dot{H}^s$ is defined by
$$\|u\|_{\dot{H}^s} =\||\partial_{x}|^s u\|_{L^2_x}.$$

Let $U(t)$ denote the solution operator to the linear Schr\"odinger equation
$$iu_t + \Delta u = 0,\; \; x \in \mathbb R,$$
that is
$$U(t)u_{0}(x)=\int e^{2\pi i\xi x-(2\pi \xi)^{2}it}\widehat{u_0}(\xi)d\xi.$$
We denote by $X^{s,b} = X^{s,b}(\mathbb R \times \mathbb R)$ the completion of
${\mathcal S}(\mathbb R \times \mathbb R)$ with respect to the following norm
$$\|u\|_{X^{s,b}}=
\|U(-t)u\|_{H_{x}^{s}H_{t}^{b}}
=\|\langle\xi\rangle^{s}\langle\tau +4\pi
^{2}\xi^{2}\rangle^{b}\tilde{u}(\xi,\tau)\|_{L_{\tau}^2L^{2}_{\xi}},$$
where $\tilde{u}(\xi,\tau)$ is the space-time Fourier Transform
$$\tilde{u}(\xi,\tau)=
\int_{\mathbb R} \int_{\mathbb R}
e^{-2 \pi i(\xi \cdot x +\tau t)}u(x,t)dxdt.$$
Furthermore for a given time interval $J$, we define
$$\|f\|_{X^{s,b}(J)}=\inf_{g  = f \textit{ on } J}  \|g\|_{X^{s,b}}.$$
Often we will drop $J$.

%

\subsection {Some known estimates} Now we recall a few known
estimates
that we shall need. First, we state the following Strichartz estimate \cite{gv}, \cite{kt}. We recall that
a pair of exponents $(q,r)$ is called admissible  in
    $\mathbb R$ if
$$ \frac{2}{q}+\frac{1}{r} = \frac{1}{2}, \ \ \ 2 \leq q,r \leq \infty.$$
\begin{prop} \label{stric}
Let $(q,r)$ and $(\tilde{q},\tilde{r})$ be
any two admissible pairs. Suppose that $u$ is a solution to
\begin{align}
&iu_{t}+ \Delta u -G(x,t)=0, \; \; x \in J \times \mathbb R ,\nonumber\\
&u(x,0)=u_{0}(x).\nonumber
\end{align} Then we have the estimate \begin{equation}\label{S}
\|u\|_{L^q_tL^r_x (J\times \mathbb R)} \lesssim
\|u_0\|_{L^2(\mathbb R)}+\|G\|_{L^{\tilde q^{\prime}}_tL^{\tilde r^{\prime}}_x
(J\times\mathbb R)}
\end{equation}
with the prime exponents denoting H\"older dual exponents.
\end{prop}

Since
$$\widehat{U(t)u_{0}}(\xi)=e^{-(2\pi \xi)^{2}it}\widehat{u_{0}}(\xi)$$
we have that
$$\|U(t)u_{0}\|_{L_{t}^{\infty}L_{x}^{2}} \lesssim \|u_{0}\|_{L^{2}}.$$
Hence
\begin{equation} \label{inf2}
\|u\|_{L_{t}^{\infty}L_{x}^{2}}
=\|U(t)U(-t)u\|_{L_{t}^{\infty}L_{x}^{2}}
\lesssim \|U(-t)u\|_{L^{\infty}_t L_{x}^{2}} \lesssim
\|u\|_{X^{0,1/2+}},
\end{equation}
where in the last inequality we applied the definition of the
$X^{s,b}$ spaces, the basic estimate $\|u\|_{L^{\infty}} \leq
\|\hat{u}\|_{L^{1}}$, and the Cauchy-Schwartz inequality.
The estimate \eqref{inf2} combined with the Sobolev embedding
theorem implies that
\begin{equation} \label{inf}
\|u\|_{L_{t}^{\infty}L_{x}^{\infty}} \lesssim
\|u\|_{X^{1/2+,1/2+}}.
\end{equation}
The Strichartz estimate gives us
\begin{equation} \label{L6}
\|u\|_{L_{t}^{6}L_{x}^{6}} \lesssim \|u\|_{X^{0,1/2+}}.
\end{equation}
If we interpolate equations
\eqref{inf} and \eqref{L6} we get
\begin{equation} \label{int1}
\|u\|_{L_{t}^{p}L_{x}^{p}} \lesssim \|u\|_{X^{\alpha_1(p),1/2+}},
\end{equation}
with $\alpha_1(p) = (\frac{1}{2}-\frac{3}{p})+$ and $6 \leq p \leq \infty$.


\section{The I-method and the local well-posedness for the $I$-system} \label{sec-I}

\subsection{The I-operator and the hierarchy of energies}


Let us define the operator $I$. For $s<1$ and a parameter $N
>>1$ let $m(\xi)$ be the following smooth monotone multiplier:
\[m(\xi):= \left\{\begin{array}{ll}
1 & \mbox{if $|\xi|<N$,}\\
(\frac{|\xi|}{N})^{s-1} & \mbox{if $|\xi|>2N$.}
\end{array}
\right.\] We define the multiplier operator $I:H^{s} \rightarrow
H^{1}$ by
$$\widehat{Iu}(\xi)=m(\xi)\hat{u}(\xi).$$
The operator $I$ is smoothing of order $1-s$ and we have that:
\begin{equation}
\|u\|_{X^{s_{0},b_{0}}} \lesssim \|Iu\|_{X^{s_{0}+1-s,b_{0}}}
\lesssim N^{1-s}\|u\|_{X^{s_{0},b_{0}}},
\end{equation}
for any $s_{0},b_{0}\in {\mathbb R}$.
\\
\\
We set
\begin{equation}\label{first}E^{1}(u)=E(Iu),\end{equation} where $$E(u)(t)=\frac{1}{2}\int
|u_{x}(t)|^{2}dx+\frac{1}{6}\int |u(t)|^{6}dx=E(u_{0}).$$ We call
$E^{1}(u)$ the first modified energy\footnote{One can actually informally define a hierarchy of modified
energies for different nonlinear dispersive equations, see \cite{ckstt1, ckstt2, ckstt4} .}.
Since we base our approach on the
analysis of a second modified energy, we collect
some facts concerning the calculus of multilinear forms used to
define the hierarchy, see, for example \cite{tz}.

If $n\geq2$ is an even integer we define a spatial multiplier of
order $n$ to be the function $M_{n}(\xi_{1},\xi_{2},\ldots ,\xi_{n})$ on
$\Gamma _{n}=\{(\xi_{1},\xi_{2},\ldots ,\xi_{n} ) \in
\{\mathbb R^{n}: \xi_{1}+\xi_{2}+ \ldots +\xi_{n}=0\}$ which
we endow with the standard measure $\delta(\xi_{1}+\xi_{2}+ \ldots
+\xi_{n})$. If $M_{n}$ is a multiplier of order $n$, $1 \leq j \leq
n$ is an index and $l \geq 1$ is an even integer,
the elongation $X_{j}^{l}(M_{n})$ of $M_{n}$ is defined to be the multiplier of
order $n+l$ given by
$$X_{j}^{l}(M_{n})(\xi_{1},\xi_{2},\ldots ,\xi_{n+l})=
M_{n}(\xi_{1},\ldots,\xi_{j-1},\xi_{j}+\ldots+\xi_{j+l},\xi_{j+l+1}, \ldots,
\xi_{n+l}).$$
Also if $M_{n}$ is a multiplier of order $n$
and $f_{1},f_{2},...,f_{n}$ are functions on $\mathbb R$ we define
$$\Lambda_{n}(M_{n};f_{1},f_{2},...,f_{n})=\int_{\Gamma_{n}}M_{n}(\xi_{1},\xi_{2},\ldots ,\xi_{n})
\prod_{i=1}^{n}\hat f_{j}(\xi_{j}),$$
and we adopt the notation
$\Lambda_{n}(M_{n};f)=\Lambda_{n}(M_{n};f,\bar f,...,f,\bar f)$.
Observe that  $\Lambda_{n}(M_{n};f)$ is invariant under
permutations of the even $\xi_{j}$ indices, or of the odd $\xi_{j}$
indices.

If $f$ is a solution of \eqref{ivp1} the following differentiation
law holds for the multilinear forms $\Lambda_{n}(M_{n};f)$:
\begin{equation}\label{diff}
\partial_{t}\Lambda_{n}(M_{n})=i\Lambda_{n}(M_{n}\sum_{j=1}^{n}(-1)^{j}\xi_{j}^{2})-i\Lambda_{n+4}
(\sum_{j=1}^{n}(-1)^{j}X_{j}^{4}(M_{n})).
\end{equation}
Observe that in this notation the first modified energy \eqref{first}
can be written as follows:
$$E^1(u)=
\frac{1}{2}\int_{\mathbb R}|\partial_{x}Iu|^{2}dx+\frac{1}{6}\int_{\mathbb R}|Iu|^{6}dx
=-\frac{1}{2}\Lambda_{2}(m_{1}\xi_{1}m_{2}\xi_{2})+\frac{1}{6}
\Lambda_{6}(m_{1}...m_{6})$$ where $m_{j}=m(\xi_{j})$.

Now we define the second modified energy
$$E^{2}(u)=-\frac{1}{2}\Lambda_{2}(m_{1}\xi_{1}m_{2}\xi_{2})
+\frac{1}{6}\Lambda_{6}(M_{6}(\xi_{1},\xi_{2},...,\xi_{6})),$$
where $M_{6}(k_{1},k_{2},...,k_{6})$ is the following multiplier:
\begin{equation}\label{M6}M_{6}(\xi_{1},\xi_{2},...,\xi_{6})=-\frac{1}{3}
\frac{m_{1}^{2}\xi_{1}^{2}-m_{2}^{2}\xi_{2}^{2}+
m_{3}^{2}\xi_{3}^{2}-m_{4}^{2}\xi_{4}^{2}+m_{5}^{2}\xi_{5}^{2}
-m_{6}^{2}\xi_{6}^{2}}{\xi_{1}^{2}-\xi_{2}^{2}+\xi_{3}^{2}-\xi_{4}^{2}+\xi_{5}^{2}-\xi_{6}^{2}}.\end{equation}
We remark that the zero set of the denominator corresponds to the
resonant set for a six-waves interaction. Also $M_6$
contains more ``cancellations" than the multiplier $m_1...m_6$
that appears in $E^1$.

The differentiation rules \eqref{diff} together with the
fundamental theorem of calculus implies the following Lemma, which
will be used to prove that $E^2$ is almost conserved.

\begin{lem}\label{fundcal}
Let $u$ be an $H^{1}$ solution to $(\ref{ivp1})$. Then for any
$T\in \mathbb R$ and $\delta >0$ we have
\begin{equation} \label{increment}
E^{2}(u(T+\delta))-E^{2}(u(T))=\int_{T}^{T+\delta}\Lambda_{10}(M_{10};u(t))dt,
\end{equation}
with $M_{10}=-\frac{i}{5!6!} \sum \{
M_{6}(\xi_{abcde},\xi_{f},\xi_{g},\xi_{h},\xi_{i},\xi_{j})
-M_{6}(\xi_{a},\xi_{bcdef},\xi_{g},\xi_{h},\xi_{i},\xi_{j})+ \\ \\
M_{6}(\xi_{a},\xi_{b},\xi_{cdefg},\xi_{h},\xi_{i},\xi_{j})-M_{6}(\xi_{a},\xi_{b},\xi_{c},\xi_{defgh},\xi_{i},\xi_{j})+\\ \\
M_{6}(\xi_{a},\xi_{b},\xi_{c},\xi_{d},\xi_{efghj},\xi_{j})-M_{6}(\xi_{a},\xi_{b},\xi_{c},\xi_{d},\xi_{e},\xi_{fghij})\}$,\\ \\
where the summation runs over all permutations
$\{a,c,e,g,i\}=\{1,3,5,7,9\}$ and $\{b,d,f,h,j\}
\\=\{2,4,6,8,10\}$. Furthermore if $|\xi_{j}|\ll N$ for all $j$ then the multiplier $M_{10}$ vanishes.
\end{lem}

As it was observed in \cite{tz} one has

\begin{prop}\label{Multb}
The multiplier $M_6$ defined in $(\ref{M6})$ is bounded on its
domain of definition.
\end{prop}

\subsection{Modified Local Well-Posedness}

In this subsection we shall prove a local well-posedness result for the
modified solution $Iu$ and some {\em a priori} estimates for it.

Let $J=[t_{0},t_{1}]$ be an interval of time. We denote by $Z_I(J)$ the
following space:
$$Z_I(J) = S_{I}(J) \cap X_{I}^{1,b}(J)$$ where $b=1/2+$ and
\begin{align*}
& S_{I}(J) = \{ f \; | \; \sup_{(q,r) \ admissible}
\|\langle \partial_{x}\rangle
If\|_{L^q_tL^r_x(J \times {\mathbb R})} < \infty \},\\
& X_{I}^{1,b}(J) = \{ f | \; \; \|If\|_{X^{1,b}(J \times {\mathbb R})} < \infty \}.
\end{align*}
\begin{prop}\label{lwp0}
Let $ s > 0$.  and consider the IVP
\begin{align}
&iIu_{t}+ \Delta Iu -I (|u|^{4}u) =0, \label{Iivp1}\\
&Iu(x,t_{0})=Iu_{0}(x)\in H^{1}({\mathbb R}),  \;\;\;\;  x \in
{\mathbb R}, \; t\in {\mathbb R}. \label{Ibc1}
\end{align}
Then for any $u_0 \in H^s$ there exists a time interval $J=[t_0,t_0+\delta], \, \delta=\delta(\|Iu_0\|_{H^1})$ and there exists a unique $u\in Z_I(J),  $ solution to \eqref{Iivp1} and \eqref{Ibc1}. Moreover there is continuity with respect to the initial data.
\end{prop}
\begin{proof}
The proof of this proposition proceeds by the usual fixed point method (see  \cite{ckstt4} ) on the space
$ Z_I(J)$. Since the estimates are very similar to the ones we provide in the proof of Proposition \ref{lwp} below,  in particular  \eqref{str-fin} and \eqref{Xsb-fin}, we omit the details. We only remark that the smallness of $\delta$ will be used once one compares \eqref{str-fin}  and \eqref{Xsb-fin} with the fact that  for any $v\in  Z_I(J)$
\[ \| Iv \|_{L^6_tL^6_x(J \times {\mathbb R})}\leq \delta^{1/6}\|Iv\|_{H^1}.\]

\end{proof}

\begin{prop}\label{lwp}
Let $ s > 0$. If $u$ is a solution to the IVP
\eqref{Iivp1}-\eqref{Ibc1} on the interval $J = [t_{0},t_{1}]$,
which satisfies the following a priori bound
$$ \|Iu\|_{L^6_t L^6_x(J \times {\mathbb R})}^6 < \mu, $$
where $\mu$ is a small universal constant,
then
\begin{equation}\label{bound1}
\|u\|_{Z_I(J)} \lesssim \|Iu_{0}\|_{H^1}.
\end{equation}
\end{prop}

\begin{proof}

We start by obtaining a control of the Strichartz norms. Applying
$\langle \partial_x \rangle $ to \eqref{Iivp1} and using the
Strichartz estimate in Proposition \ref{stric}, for any pair of
admissible exponents $(q,r)$ we obtain
\begin{equation} \label{str-duh}
\|\langle \partial_x \rangle Iu\|_{L^q_tL^r_x}
\lesssim \|Iu_0\|_{H^1_x}
+ \|\langle \partial_x \rangle I(|u|^4 u)\|_{L^{\frac{6}{5}}_t L^{\frac{6}{5}}_x}.
\end{equation}
Now we notice that the multiplier $\langle \partial_x \rangle I$
has a symbol which is increasing as a function of $|\xi|$ for any
$s\geq 0$. Using this fact one can modify the proof of the
Leibnitz rule for fractional derivatives and prove its validity
for  $\langle \partial_x \rangle I$. See also Principle $A. 5$ in the appendix of \cite{tt}. This remark combined with
\eqref{str-duh} implies that
\begin{align}
\|\langle \partial_x \rangle Iu\|_{L^q_tL^r_x}
& \lesssim \|Iu_0\|_{H^1_x}
+ \|\langle \partial_x \rangle Iu\|_{L^6_t L^6_x}
\|u^4\|_{L^{\frac{3}{2}}_t L^{\frac{3}{2}}_x} \nonumber \\
&  \lesssim \|Iu_0\|_{H^1_x}
+ \|u\|_{Z_{I}(J)}
\|u\|^4_{L^6_t L^6_x}, \label{str}
\end{align}
where to obtain \eqref{str} we used H\"{o}lder's inequality and the definition of
$Z_I(J)$.

In order to obtain an upper bound on $\|u\|_{L^6_tL^6_x}$ we perform a
Littlewood-Paley decomposition along the following lines. We note that a similar approach was used in \cite{ckstt4}, Lemma 3.1. We write
\begin{equation} \label{lp}
u = u_{N_0} + \sum_{j=1}^{\infty} u_{N_j},
\end{equation}
where $u_{N_0}$ has spatial frequency support for $\langle \xi
\rangle \leq N$, while $u_{N_j}$ is such that its spatial Fourier
transform is supported for $\langle \xi \rangle \sim N_j =2^{h_j}$
with $h_j \gtrsim \log N$ and $j = 1,2,3,...$. Let $\epsilon > 0$.
The triangle inequality applied on \eqref{lp} gives
\begin{align}
\|u\|_{L^6_t L^6_x}
& \leq \|u_{N_0}\|_{L^6_tL^6_x} + \sum_{j=1}^{\infty} \|u_{N_j}\|_{L^6_tL^6_x} \nonumber \\
& = \|Iu_{N_0}\|_{L^6_tL^6_x} + \sum_{j=1}^{\infty}
\|u_{N_j}\|_{L^6_tL^6_x}^{\epsilon}
\|u_{N_j}\|_{L^6_tL^6_x}^{1-\epsilon}. \label{lp-tri}
\end{align}
On the other hand, by using the definition of the operator $I$,
the definition of the $u_{N_j}$'s  and
the Marcinkiewicz multiplier theorem
we observe the following:
\begin{align}
& \|Iu_{N_0}\|_{L^6_tL^6_x} \lesssim \|Iu\|_{L^6_tL^6_x}, \label{inI} \\
& \|u_{N_j}\|_{L^6_tL^6_x} \lesssim N_j^{1-s} N^{s-1}
\|Iu_{N_j}\|_{L^6_tL^6_x},\;\; j=1,2,3,... \label{Iuj}\\
& \|u_{N_j}\|_{L^6_tL^6_x} \lesssim N_j^{-s} N^{s-1} \|\langle
\partial_x \rangle Iu_{N_j}\|_{L^6_tL^6_x},\;\; j=1,2,3,...
\label{DIuj}
\end{align}
Now we use the estimates \eqref{inI} - \eqref{DIuj} to obtain the following
upper bound on \eqref{lp-tri}
$$
\|u\|_{L^6_t L^6_x} \lesssim  \|Iu\|_{L^6_tL^6_x} +
\sum_{j=1}^{\infty} N_j^{-s+\epsilon} N^{s-1}
\|Iu_{N_j}\|_{L^6_tL^6_x}^{\epsilon} \|\langle \partial_x \rangle
Iu_{N_j}\|_{L^6_tL^6_x}^{1-\epsilon},
$$
which after noticing that $N^{s-1} \leq 1$ implies
\begin{equation} \label{LPready}
\|u\|_{L^6_t L^6_x} \lesssim  \|Iu\|_{L^6_tL^6_x} +
\sum_{j=1}^{\infty} N_j^{-s+\epsilon}
\|Iu_{N_j}\|_{L^6_tL^6_x}^{\epsilon} \|\langle \partial_x \rangle
Iu_{N_j}\|_{L^6_tL^6_x}^{1-\epsilon}.
\end{equation}
Now using the cheap Littlewood-Paley inequality
$$\sup_{j}\|u_{N_{j}}\|_{L^p} \lesssim \|u\|_{L^p}$$
for $1 \leq p \leq \infty$, and summing \eqref{LPready}, we have that for any $s > \epsilon$
\begin{equation} \label{L66lp}
\|u\|_{L^6_t L^6_x}
\lesssim  \|Iu\|_{L^6_tL^6_x} +
\|Iu\|_{L^6_tL^6_x}^{\epsilon}
\|\langle \partial_x \rangle Iu\|_{L^6_tL^6_x}^{1-\epsilon},
\end{equation}
which combined with \eqref{str} implies
\begin{equation} \label{str-fin}
\|\langle \partial_x \rangle Iu\|_{L^q_tL^r_x} \lesssim
\|Iu_0\|_{H^1_x} + \|u\|_{Z_I(J)} \|Iu\|^4_{L^6_tL^6_x} +
\|u\|_{Z_I(J)}^{5-4\epsilon} \|Iu\|_{L^6_tL^6_x}^{4\epsilon}.
\end{equation}

Now we shall obtain a control of the $X^{s,b}$ norm.
We use Duhamel's formula and the theory of $X^{s,b}$ spaces (for details, see, \cite{g}, \cite{kpv}) to obtain

\begin{equation} \label{Xsb-duh}
\|\langle \partial_x \rangle Iu\|_{X^{0, \frac{1}{2}+}} \leq \|Iu_0\|_{H^1_x} + \|\langle \partial_x \rangle I(|u|^4u)\|_{X^{0, -\frac{1}{2}+}}.
\end{equation}
However interpolating between
$$ \|u\|_{L^6_tL^6_x} \lesssim \|u\|_{X^{0,\frac{1}{2}+}}$$
and
$$ \|u\|_{L^2_tL^2_x} \lesssim \|u\|_{X^{0,0}}$$
we obtain
$$ \|u\|_{L^{6-}_tL^{6-}_x} \lesssim \|u\|_{X^{0, \frac{1}{2}-}},$$
which by duality gives
$$ \|u\|_{X^{0, -\frac{1}{2} + }} \lesssim \|u\|_{L^{\frac{6}{5}+}_tL^{\frac{6}{5}+}_x}.$$
Hence \eqref{Xsb-duh} implies
\begin{equation} \label{Xsb-duh2}
\|Iu\|_{X^{1, \frac{1}{2}+}} \lesssim \|Iu_0\|_{H^1_x}
+ \|\langle \partial_x \rangle I(|u|^4 u)\|_{L^{\frac{6}{5}+}_t L^{\frac{6}{5}+}_x}.
\end{equation}
As we noticed earlier when obtaining an upper bound on the Strichartz norms,
the Leibnitz rule for fractional derivatives can be proved for
$\langle \partial_x \rangle I$. Therefore after applying the
Leibnitz rule to \eqref{Xsb-duh2} we obtain:
\begin{align}
\|Iu\|_{X^{1, \frac{1}{2}+}}
& \lesssim \|Iu_0\|_{H^1_x}
+ \|\langle \partial_x \rangle Iu\|_{L^6_t L^6_x}
\|u^4\|_{L^{\frac{3}{2}+}_t L^{\frac{3}{2}+}_x} \nonumber \\
&  \lesssim \|Iu_0\|_{H^1_x}
+ \|u\|_{Z_{I}(J)}
\|u\|_{L^{6+}_t L^{6+}_x} \|u\|_{L^6_tL^6_x}^3,\label{Xsb}
\end{align}
where to obtain \eqref{Xsb} we used H\"{o}lder's inequality and the definition of
$Z_I(J)$.
An upper bound on $\|u\|_{L^6_tL^6_x}$ is given by \eqref{L66lp}.
In order to obtain an upper bound on $\|u\|_{L^{6+}_tL^{6+}_x}$ we proceed
as follows. First, we perform a dyadic decomposition and write
$u$ as in \eqref{lp}. The triangle inequality applied on \eqref{lp} gives
\begin{align}
\|u\|_{L^{6+}_t L^{6+}_x}
& \lesssim \|u_{N_0}\|_{L^{6+}_tL^{6+}_x} + \sum_{j=1}^{\infty} \|u_{N_j}\|_{L^{6+}_tL^{6+}_x} \nonumber \\
& = \|Iu_{N_0}\|_{L^{6+}_tL^{6+}_x} + \sum_{j=1}^{\infty}
N_j^{\delta - s} N^{s-1} \| \langle \partial_x \rangle^{1-\delta}
Iu_{N_j}\|_{L^{6+}_tL^{6+}_x}, \label{lpI}
\end{align}
where in order to obtain \eqref{lpI} we used the following
inequality which can be verified
by recalling the definition of the operator $I$, the definition of the $u_{N_j}$'
and the Marcinkiewicz multiplier theorem:
$$ \|u_{N_j}\|_{L^{6+}_tL^{6+}_x} \lesssim N_j^{\delta -s} N^{s-1}
\|\langle \partial_x \rangle^{1-\delta}
Iu_{N_j}\|_{L^{6+}_tL^{6+}_x},\;\; j=1,2,3,...$$ In order to
further bound the right hand side of \eqref{lpI}, first we notice
that $N^{s-1} \leq 1$. Then we apply the cheap Littlewood-Paley inequality as before
and sum \eqref{lpI} to obtain for
any $s > \delta$
\begin{equation} \label{L66lp+}
\|u\|_{L^{6+}_t L^{6+}_x}
\lesssim  \|Iu\|_{L^{6+}_tL^{6+}_x} +
\|\langle \partial_x \rangle^{1-\delta} Iu\|_{L^{6+}_tL^{6+}_x}.
\end{equation}
However
\begin{equation} \label{6+}
\|\langle \partial_x \rangle^{1 - \delta} Iu\|_{L^{6+}_t L^{6+}_x} \leq \|Iu\|_{X^{1,\frac{1}{2}+}},
\end{equation}
which follows from interpolating
$$ \|u\|_{L^6_tL^6_x} \leq \|u\|_{X^{0,\frac{1}{2}+}}$$
and
$$ \|u\|_{L^{\infty}_tL^{\infty}_x} \leq \|u\|_{X^{\frac{1}{2}+,\frac{1}{2}+}}.$$
Now we combine \eqref{L66lp+} and \eqref{6+} to obtain
\begin{equation} \label{L6+ready}
\|u\|_{L^{6+}_t L^{6+}_x}
\lesssim  \|Iu\|_{X^{1,\frac{1}{2}+}}.
\end{equation}

By applying the inequalities obtained in \eqref{L66lp} and \eqref{L6+ready} to bound the right hand side of \eqref{Xsb},
we obtain
\begin{equation} \label{Xsb-fin}
\|Iu\|_{X^{1, \frac{1}{2}+}} \lesssim \|Iu_0\|_{H^1_x} +
\|u\|_{Z_{I}(J)}^2 \|Iu\|_{L^6_tL^6_x}^3 +
\|u\|_{Z_I(J)}^{5-3\epsilon} \|Iu\|_{L^6_tL^6_x}^{3\epsilon}.
\end{equation}

The desired bound \eqref{bound1}
follows from \eqref{str-fin} and \eqref{Xsb-fin}.
\end{proof}

\subsection{${\mathbf{E^{2}(u)}}$ is a small perturbation of ${\mathbf{E^{1}(u)}}$ }
Now we shall prove that the second energy $E^2(u)$ is a small perturbation of the first
energy $E^1(u)$.

\subsection*{Decomposition remark} Our approach to prove that the second energy is
a small perturbation of the first energy as well as to prove a
decay for the increment of the second energy is based on obtaining certain
multilinear estimates in appropriate functional spaces which are
$L^2$-based. Hence, whenever we perform a Littlewood-Paley
decomposition of a function we shall assume that the Fourier
transforms of the Littlewood-Paley pieces are positive. Moreover,
we will ignore the presence of conjugates. At the end we will always keep a factor of order $N_{max}^{-\epsilon}$, where
 $N_{max}$ is the largest frequency of the different pieces, in order to perform the summations.
The details for both the Propositions \ref{e1} and \ref{e2} can be found in \cite{tz},
but we choose to present the main points of the argument to make the paper as self-contained as possible.

\begin{prop}\label{e1}
Assume that $u$ solves $(\ref{ivp1})$ with $s>1/3$. Then
$$E^{2}(u) = E^{1}(u)+O(N^{-\epsilon})\|Iu\|_{H^{1}}^{6}.$$
Moreover if $\|Iu\|_{H^{1}}=O(1)$ then
$\|\partial_{x}Iu\|_{L^{2}}^{2}\lesssim E^{2}(u)$.
\end{prop}
\begin{proof}

By the definition of the first and second modified energy we have
\begin{equation} \label{E12}
E^{2}(u)=-\frac{1}{2}\Lambda_{2}(m_{1}\xi_{1}m_{2}\xi_{2})+
\frac{1}{6} \Lambda_{6}(M_{6})=E^{1}(u)+
\frac{1}{6}\Lambda_{6}(M_{6}-\prod_{i=1}^{6}m_{i} ).
\end{equation}
Hence it suffices to prove the following pointwise in time
estimate
\begin{equation}\label{perturbation}
|\Lambda_{6}(M_{6}-\prod_{i=1}^{6}m_{i} )|(t) \lesssim O(N^{-\epsilon})
\|Iu(\cdot, t)\|_{H^{1}_x}^{6}.\end{equation} Combining the
\textbf{Decomposition Remark} with the fact that $M_6$ is bounded
(by Proposition \ref{Multb}) and that $m$ is bounded (by its
definition) it is enough to show that
\begin{equation}\nonumber \int_{\Gamma_{6}}
\prod_{j=1}^{6}\hat{u}(\xi_{j}) \lesssim \frac{1}{N^{\epsilon}}\|Iu(\cdot,
t)\|_{H^{1}_x}^{6}.\end{equation}
Towards this aim, we use again a dyadic decomposition and one can
easily check that it suffices to show the following:
\begin{equation}\label{pieceb}
\int_{\Gamma_{6}} \prod_{j=1}^{6}\widehat{u_{N_{j}}}(\xi_{j}) \lesssim
\frac{1}{N^{\epsilon}}(N_1...N_6)^{0-}\prod_{j=1}^{6}\|Iu_{N_{j}}(\cdot,
t)\|_{H^{1}_x}^{6}\end{equation} where we recall that $u_{N_{j}}$ is supported around
$\langle \xi \rangle \sim N_j$.

 We can rearrange the sizes of the frequencies so that  $N_{1}^{*} \geq N_{2}^{*} \geq N_{3}^{*}
\geq N_{4}^{*} \geq N_{5}^{*} \geq N_{6}^{*}$ and for simplicity set
$u_{N_{j}}=u_{j}$ and $u_{N_{j}^{*}}=u_{j}^{*} $.

We may assume that $N_{1}^* \gtrsim N$, otherwise
$M_{6}-\prod_{i=1}^{6}m_{i} \equiv 0$ and \eqref{perturbation}
follows trivially. We have that
$$\int_{\Gamma_{6}} \prod_{j=1}^{6}\widehat{u_j}(\xi_{j}) \lesssim
\frac{1}{(N_1^{*})^\epsilon}\int \widehat{
(N_1^{*})^{\epsilon}u_1^*}\prod_{j=2}^{6}\widehat{u_j^*} \lesssim
\frac{1}{N^{\epsilon}}(\||\nabla|^{\epsilon}u_1^*\|_{L_x^{6-}})
\prod_{j=2}^{6}\|u_j^*\|_{L_{x}^{6+}}$$ by reversing
Plancherel and applying H\"older's inequality. Moreover by Sobolev
embedding
$$\||\nabla|^{\epsilon}u_1^*\|_{L_x^{6-}}\lesssim \|u_j^*\|_{H^{1/3+}_x}$$
and thus
$$\int_{\Gamma_{6}} \prod_{j=1}^{6}\widehat{u_j}(\xi_{j}) \lesssim
\frac{1}{N^{\epsilon}}\|u_j^*\|_{H^{1/3+}_x}^{6}$$
In addition for $s>1/3$ we have
$$\|u_j^*\|_{H^{1/3+}_x} \lesssim \|Iu_j^*\|_{H^{1+\frac{1}{3}+\epsilon -s}_x} \lesssim
\|Iu_j^*\|_{H^{1}_x}.$$
Thus \eqref{pieceb} follows.
\end{proof}

{\bf Remark.} As mentioned in the introduction, we obtain global well-posedness for $s>1/3$, which is exactly the restriction on $s$
of the previous proposition. So this regularity is a threshold to the method, at least if
 one wants to use the second modified energy.

\subsection{An upper bound on the increment of ${\mathbf{E^2(u)}}$}

In Lemma \ref{fundcal} we proved that an increment of the second
energy can be expressed as
$$
E^{2}(u(T+\delta))-E^{2}(u(T))=\int_{T}^{T+\delta}\Lambda_{10}(M_{10};u(t))dt.$$
Hence in order to control the increment of the second energy we shall find an
upper bound on the $\Lambda_{10}$ form, which we do in the following
proposition.

\begin{prop}\label{e2}
For any Schwartz function u, and any $\delta \sim 1$, we have that
$$\left|\int_{0}^{\delta}\Lambda_{10}(M_{10};u(t))\right| \lesssim N^{-2+}\|Iu\|_{X^{1,1/2+}}^{10},$$ for $s>1/4.$
\end{prop}

\begin{proof}
We observe that $M_{10}$ is bounded as an elongation of the
bounded multiplier $M_6$. We perform a dyadic decomposition as in
Proposition \ref{e1}, and we borrow the same notation. Since we are integrating over
$\Gamma_{10}$, we have that $N_1^* \sim N_2^*$, and we may assume
that $N^*_1 \gtrsim N$ otherwise $M_{10} \equiv 0.$ Thus we consider $N_{1}^{*} \sim N_{2}^{*} \gtrsim N$.
\\
\\
Since
$\frac{1}{m(N_{1}^{*})m(N_{2}^{*})(N_{1}^{*}N_{2}^{*})^{1-}}
\lesssim \frac{1}{N^{2-}}$ we have
\begin{align}
|\int_{0}^{\delta}\int M_{10}\prod_{j=1}^{10}\hat{u}_{j}|
&\lesssim \frac{(N_1^*)^{0-}}{N^{2-}}
\|\langle \partial_{x} \rangle ^{1-}Iu_{1}^*\|_{L_{t}^{6}L_{x}^{6}}
\|\langle \partial_{x} \rangle^{1-}¥ Iu_{2}^*\|_{L_{t}^{6}L_{x}^{6}}
\prod_{j=3}^{10}\|u_{j}^*\|_{L_{t}^{12}L_{x}^{12}} \nonumber
\\ 
& \lesssim
\frac{(N_1^*)^{0-}}{N^{2-}}\|Iu\|_{X^{1,1/2+}}^{2}\|u\|_{X^{(\frac{1}{2}-\frac{3}{12})+,1/2+}}^{8}
\label{useinter}\\ & \lesssim
(N_1^*)^{0-}N^{-2+}\|Iu\|_{X^{1,1/2+}}^{10}\label{spaces}
\end{align} where in order to obtain \eqref{useinter} we use
\eqref{L6} and to obtain \eqref{spaces} we use the fact that for
$s>\frac{1}{4}$ the following inequality holds
$$\|u\|_{X^{(\frac{1}{2}-\frac{1}{4})+,1/2+}}=\|u\|_{X^{1/4+,1/2+}}\lesssim \|Iu\|_{X^{1,1/2+}}.$$

\end{proof}

\section{The almost Morawetz estimate} \label{sec-almor}

In this section we shall prove an interaction Morawetz estimate
for the smoothed out solution $Iu$, hence the name ``almost
Morawetz" estimate. Our estimate reads as follows:


\begin{thm}\label{main}Let $u \in \mathcal{S}$ be a solution to the NLS
\begin{equation}\label{NLSq}
iu_{t}+\Delta u=|u|^4u, \ \ \ (x,t) \in \mathbb{R} \times [0,T].
\end{equation} Then, there exists a convex function $a: \R^4 \rightarrow \R$
with $\nabla a \in L^{\infty}(\R^4)$, such that
\begin{align}\label{amorerror}
\|Iu\|_{L_{t}^{8}L_{x}^{8}([0,T] \times {\mathbb R} }^{8}
&\lesssim \sup_{[0,T]}\|Iu\|_{\dot{H}^{1}_x}\|Iu\|_{L^{2}_x}^{7}
\\ \nonumber & + \left|\int_{0}^{T}\int_{\Bbb R^4}\nabla
a(x_{1},x_{2},x_{3},x_{4}) \cdot
\{\CN_{bad},Iu_1Iu_2Iu_3Iu_4\}_{p}dx_1dx_2dx_3dx_4dt\right|,
\end{align}
where
$$\CN_{bad}=\sum_{k=1}^{4}\left(I(\CN_{k}(u))-\CN_{k}(Iu)\right)\prod_{j=1,j \ne k}^{4}Iu_{j},$$
$\CN(f) = |f|^4f$, $f_k = f(x_k, t)$, $\CN_k(f) = \CN(f_k)$, and
$\{\cdot\}_{p}$ is the momentum bracket defined by
$$\{f,g\}_p=\Re(f\overline{\nabla g} -
g\overline{\nabla f}).$$Moreover in a time interval
$J=[t_{0},t_{1}]$ where the solution $u$ belongs to the space
$Z_I(J)$ we have that
\begin{equation}\label{errorcontrol}
\left|\int_{t_0}^{t_1}\int_{\Bbb R^4}\nabla
a(x_{1},x_{2},x_{3},x_{4}) \cdot
\{\CN_{bad},Iu_1Iu_2Iu_3Iu_4\}_{p}dx_1dx_2dx_3dx_4dt\right|
\lesssim \frac{1}{N^-}\|u\|^{12}_{Z_I(J)}.
\end{equation}
\end{thm}

\begin{rem}
In all of our arguments we will work with smooth (Schwartz)
solutions. This will simplify the calculations and will enable us
to justify the steps in the subsequent proofs. Then the local
well-posedness theory and the perturbation theory that have been
established (see, for example, \cite{cw1}) for this problem can be
applied to approximate the $H^{s}$ solutions by smooth solutions.
\end{rem}

We prove the above almost Morawetz estimate inspired by the idea
of the proof of the interaction Morawetz estimate for the
defocusing nonlinear cubic Schr\"odinger equation on ${\mathbb
R}^3$, \cite{ckstt4} and the one  that recently appeared in \cite{chvz06}.
However we establish a Morawetz estimate for
the almost solution, i.e for $Iu$ itself, which is the main novelty
of our approach. In order to make our presentation complete, first
we recall a general approach of obtaining interaction Morawetz
estimates \cite{Inotes}.
Then we present our derivation of the interaction
Morawetz estimate for the almost solution $Iu$.


\subsection{Towards interaction Morawetz estimates.}
In this subsection we follow \cite{Inotes}. Note that here we work
in general dimension $d$.

Let $u \in \mathcal{S}$ be a solution to the NLS
\begin{equation}\label{NLS}
iu_{t}+\Delta u=\CN(u), \ \ \ (x,t) \in \mathbb R^d \times [0,T].
\end{equation}
 We say that $\CN$
corresponds to a defocusing potential $G$ (meaning $G$ positive)
if
$$\{\CN,u\}_{p}^{j}=\partial_{j}G.$$ Here we are denoting by
$\{f,g\}_p$ the vector whose components are given by
$$\{f,g\}^j_p=\Re(f\overline{\partial_jg} -
g\overline{\partial_jf}).$$ For example, in the case when
$\CN(u)=|u|^4u$ we have that $\{\CN,u\}_{p}^{j}=-\partial_{j}G$,
where $G=\frac{2}{3}|u|^{6}$.

Now let us define the momentum density via
\begin{equation} \label{defT}
T_{0j}=2\Im (\bar{u}\partial_{j}u)
\end{equation}
and the linearized momentum current
\begin{equation} \label{defL}
L_{jk}=-\partial_{j}\partial_{k}(|u|^2)+4\Re (\overline{\partial_{j}u}\partial_{k}u).
\end{equation}
An easy computation shows that
\begin{equation} \label{derT}
\partial_{t}T_0j+\partial_{k}L_{jk}=2\{\CN,u\}_{p}^{j}\; ,
\end{equation}
where we adopt Einstein's summation convention. Thus by
integrating in space we have that in the case when $\CN$
corresponds to a potential, then the total momentum
 is conserved in time,
$$\int_{\mathbb{R}}T_{0j}(x,t)dx=C.$$
Finally, if $a:\R^d \rightarrow \R$ is convex then we define the
Morawetz action associated to $u$ by the formula
\begin{equation}
\label{Mat}M_a(t)= 2 \int_{\Bbb R^d}\nabla a(x) \cdot
\Im(\overline{u}(x)\nabla u(x))dx.
\end{equation}

We now recall a classical result. The first step in the proof of
the estimate \eqref{amorerror} is obtained by a slight
modification of the argument in the following Proposition, in the
case when the forcing term $\mathcal{N}$ does not correspond to a
defocusing potential. We will state this result in the form of a
corollary.

\begin{prop} \label{clasprop}
Let $a:\R^d \rightarrow \R$ be a convex function and $u$ be a
smooth solution to equation \eqref{NLS} on ${\mathbb R}^{d}\times
[0,T]$ with a defocusing potential $G$. Then, the following
inequality holds
\begin{equation}\label{Mor}
\int_{0}^{T}\int_{\Bbb R^d}(-\Delta \Delta a)|u(x,t)|^2dxdt
\lesssim \sup_{t \in [0,T]}|M_{a}(t)|.
\end{equation}
\end{prop}

\begin{proof} Without loss of generality, we can assume that $a$
is smooth. Then, a standard approximation argument concludes the
proof for the general case.

 According to \eqref{defT} the Morawetz
action can be written as
$$M_{a}(t)=\int_{\Bbb R^d} \partial_{j}a \; T_{0j}.$$
Then thanks to \eqref{derT},
\begin{align}
\partial_{t}M_a(t) & =\int_{\Bbb R^d} \partial_{j}a \; \partial_{t} T_{0j} \\
& =\int_{\Bbb R^d} \partial_{j}a\left( -\partial_{k}L_{jk}+2\{\CN,u\}_{p}^{j}\right) \nonumber \\
& = \int_{\Bbb R^d} \partial_{j}a\left( -\partial_{k}L_{jk}-2\partial_{j}G\right)\nonumber\\
& =\int_{\Bbb R^d}(\partial_{j}\partial_{k}a)L_{jk}dx +2\int_{\Bbb R^d}\Delta a \; G \; dx, \label{dermor}
\end{align}
where in the last equality we used integration by parts.
Now \eqref{dermor} combined with the definition  of $L_{jk}$  \eqref{defL} implies
\begin{align*}
\partial_{t}M_a(t) & =
\int_{\mathbb R^d}(-\partial_{j}\partial_{k}a)\partial_{j}\partial_{k}(|u|^2)\; dx
+ 4\int_{\mathbb R^d}(\partial_{j}\partial_{k}a)\Re \left( \overline{\partial_{j} u} \partial_{k}u\right)dx
+2\int_{\mathbb R^d}\Delta a \; G \; dx\\
& =-\int_{\Bbb R^d}(\Delta \Delta a) |u|^{2}dx
+2\int_{\Bbb R^d}\Delta a \; G \; dx
+ 4\int_{\mathbb R^d}(\partial_{j}\partial_{k}a)\Re
\left( \overline{\partial_{j} u}\partial_{k}u\right)dx.
\end{align*}
Since $\partial_{j}\partial_{k}a$ is weakly convex we have that
$$4(\partial_{j}\partial_{k}a)\Re \left( \overline{\partial_{j} u}\partial_{k}u\right) \geq 0$$
and the trace of the Hessian of $\partial_{j}\partial_{k}a$, which is $\Delta a$, is positive. Thus
$$-\int_{\Bbb R^d}(\Delta \Delta a) |u|^{2}dx \leq \partial_{t}M_a(t).$$
Hence by the fundamental theorem of calculus we have that
\begin{equation}
\int_{0}^{T}\int_{\Bbb R^d}(-\Delta \Delta a)|u(x,t)|^2dxdt
\lesssim \sup_{t \in [0,T]}|M_{a}(t)|. \label{mora}
\end{equation}
\end{proof}

In the case of a solution to an equation with a nonlinearity which
is not associated to a defocusing potential, we immediately obtain
the following corollary.

\begin{cor} \label{clasprop2}
Let $a:\R^d \rightarrow \R$ be convex and $v$ be a smooth solution
to the equation
\begin{equation}\label{NLS2} iv_{t}+\Delta v=\widetilde{\CN}, \ \
\ (x,t) \in \mathbb R^d \times [0,T].
\end{equation}Then, the following inequality
holds
\begin{equation}\label{Mor}
\int_{0}^{T}\int_{\Bbb R^n}(-\Delta \Delta a)|v(x,t)|^2dxdt +2
\int_0^T\int_{\R^d} \nabla a \cdot \{\widetilde{\CN},v\}_pdxdt
\lesssim \sup_{t \in [0,T]}|M_{a}(t)|,
\end{equation} where $M_{a}(t)$ is the Morawetz action
corresponding to $v$.
\end{cor}

\vspace{2mm}

Now we are ready to derive the main inequality in obtaining
interaction Morawetz estimates. Let $u_{k},$ $k=1,2,$ be solutions
to \eqref{NLS} with nonlinearity $\widetilde{\CN}_k$ in
$d_{k}-$spatial dimensions. Assume $\widetilde{\CN}_k$ has a
defocusing potential $G_k.$ Define the tensor product
$u:=(u_{1}\otimes u_{2})(x,t)$ for $x$ in $\mathbb R^{d_{1}}
\times \R^{d_2} $ by the formula
$$(u_{1}\otimes u_{2})(x,t)=u_{1}(x_{1},t) \; u_{2}(x_{2},t).$$
Then, it can be easily verified that $u_1 \otimes u_2$ solves
(\ref{NLS}) with forcing term $\CN=\widetilde{\CN}_{1}\otimes
u_2+\widetilde{\CN}_{2}\otimes u_1$. Since
$$\{\widetilde{\CN}_{1}\otimes u_2+\widetilde{\CN}_{2}\otimes u_1,u_1 \otimes u_2\}_{p}=\left(\{\widetilde{\CN}_1,u_1\}_{p}
\otimes |u_2|^2, \{\widetilde{\CN}_2,u_2\}_{p}\otimes
|u_1|^2\right)$$ we have the important fact that the tensor
product of defocusing semilinear Sch\"odinger equations is also
defocusing in the sense that
$$\{\widetilde{\CN}_{1}\otimes u_2+\widetilde{\CN}_{2}\otimes u_1,u_1 \otimes u_2\}_{p}=-\nabla G,$$
where  $\nabla= (\nabla_{x_1},\nabla_{x_2})$ and $G=G_{1}\otimes
|u_2|^2+G_{2}\otimes |u_1|^2$, thus $G \geq 0$. Since $u_1 \otimes
u_2$ solves (\ref{NLS}) and obeys momentum conservation with a
defocusing potential, we can apply Proposition \ref{clasprop} to
obtain
\begin{equation}\label{interact}
\int_{0}^{T}\int_{\Bbb R^{d_1} \times \Bbb R^{d_2}}(-\Delta \Delta
a)|u_1\otimes u_2|^2(x,t)dxdt \lesssim
\sup_{t \in [0,T]}|M_{a}^{\otimes_{2}}(t)|,
\end{equation}
where $\Delta= \Delta_{x_1}+\Delta_{x_2}$ is the $d_{1}+d_{2}$
Laplacian, $a$ is any real-valued convex function on $\R^{d_1}
\times \R^{d_2},$ and $M_{a}^{\otimes_{2}}(t)$ is the Morawetz
action
 that corresponds to $u_1\otimes u_2$.

Clearly, this argument can be generalized to an arbitrary number
of solutions $u_k$ to \eqref{NLS} with nonlinearity
$\widetilde{\CN}_k$ with a defocusing potential $G_k.$ Indeed one
obtains

\begin{equation}\label{interact2}
\int_{0}^{T}\int_{\Bbb R^{d_1} \times \ldots \times  \Bbb
R^{d_k}}(-\Delta \Delta a)|v|^2(x,t)dxdt \lesssim
\sup_{t \in [0,T]}|M_{a}^{\otimes_{k}}(t)|,
\end{equation}
with $x :=(x_1,\ldots, x_k), v(x):=\bigotimes_{i=1}^{k} u_i(x_i),
\widetilde{\CN}:=\sum_{i=1}^k \widetilde{\CN_i} \bigotimes_{j\neq
i} u_j ,$ and $M_{a}^{\otimes_{k}}(t)$ the Morawetz action
corresponding to $v$.

Moreover, in the case when $\widetilde{\CN}_k$ does not have a
defocusing potential, then according to Corollary \ref{clasprop2}
we get,
\begin{equation}\label{interact3}
\int_{0}^{T}\int_{\Bbb R^{d_1} \times \ldots \times  \Bbb
R^{d_k}}(-\Delta \Delta a)|v|^2(x,t)dxdt +2 \int_0^T\int_{\Bbb
R^{d_1} \times \ldots \times  \Bbb R^{d_k}} \nabla a \cdot
\{\widetilde{\CN},v\}_pdxdt \lesssim
\sup_{t \in [0,T]}|M_{a}^{\otimes_{k}}(t)|.
\end{equation}

\noindent We will use \eqref{interact3} in the proof of Theorem
\ref{main}.

\subsection{Interaction Morawetz estimates.}
Estimate \eqref{interact2} turns out to be very useful,
as a careful choice of the function $a$ allows one to obtain
bounds on a particular Lebesgue norm of a solution to equation
\eqref{NLS} with defocusing potential. For example, let $u$ be a
solution to \eqref{NLS} with a defocusing potential, in dimension
$d=3$. Let $k=2$, $d_1=d_2=3$, $u_1=u_2=u$, and choose
$$a(x)=a(x_1,x_2)=|x_1-x_2|.$$ Then an easy calculation shows that
$-\Delta \Delta a=C\delta(|x_1-x_2|)$, and equation
\eqref{interact2} gives
$$\int_{0}^{T}\int_{\Bbb R^3}|u(x,t)|^{4}dx \lesssim
\sup_{t \in [0,T]}|M_{a}^{\otimes_{2}}(t)|.$$
It can be shown using Hardy's inequality (for details see
\cite{ckstt4}) that  when $d=3$
$$\sup_{t \in [0,T]}|M_{a}^{\otimes_{2}}(t)|
\lesssim \sup_{t \in [0,T]}\|u(t)\|_{\dot{H}^{\frac{1}{2}}}^2$$
and thus
$$\|u(x,t)\|^{4}_{L^4_{t \in [0,T]}L^4_x} \lesssim \sup_{t \in [0,T]} \|u(t)\|_{\dot{H}^{\frac{1}{2}}}^2$$
which is the interaction Morawetz estimate that appears in
\cite{ckstt4}.

Analogously, in the case $d=1$, let $u$ be a solution to
\eqref{NLS} with a defocusing potential, and let $k=4$,
$d_1=\ldots=d_4=1,$ and $u_1=\ldots=u_4=u.$ Then \eqref{interact2}
reads

\begin{align}\label{holm}
\int_{0}^{T}\int_{\Bbb R^4}(-\Delta \Delta
a)\prod_{i=1}^{4}|u(x_{i},t)|^2dx_{1}dx_2dx_3dx_4dt  \lesssim
\sup_{t \in [0,T]}|M_{a}^{\otimes_{4}}(t)|.
\end{align}

In order to proceed as in the  case $d=3$  and obtain a bound on a
Lebesgue norm of $u$, we need to choose an appropriate function
$a$, and get an upper bound on the right hand side. An elegant
idea to obtain the desired bound can be found in \cite{chvz06}.
Precisely, one performs an orthonormal change of variables $z=Ax$
with $A$ orthonormal matrix, and writes \eqref{holm} with respect
to the $z$-variable. Notice that $\Delta_{z}=\Delta_{x}$ and also
the orthonormal change of variables leaves invariant the inner
product which appears in the Morawetz action on the right hand
side of \eqref{holm}. Choosing the convex function $a(z)$ to be
$a(z)=(z_{2}^{2}+z_{3}^{2}+z_{4}^{2})^{1/2}$ (hence $-\Delta_{z}
\Delta_{z} a(z)=4\pi \delta(z_{2},z_{3},z_{4})$), it is possible
to estimate quickly the right hand side, and then going back to
the $x$-variable one obtains the following estimate
\begin{equation} \label{intMor1d}
\|u\|_{L_{t \in [0,T]}^{8}L_{x}^{8}}^{8} \lesssim
\sup_{t \in [0,T]}\|u\|_{\dot{H}^{1}}\|u\|_{L^{2}}^{7}.
\end{equation}
For the details we refer the reader to \cite{chvz06}. We will use
the same orthonormal transformation and the same choice of
function $a$ in the proof of our almost Morawetz estimate.

\subsection{Almost Morawetz estimate. Proof of Theorem \ref{main}.}


Recall that $u \in \mathcal{S}$ is a solution to the NLS
\begin{equation}\label{NLSproof}
iu_{t}+\Delta u=\CN(u), \ \ \ (x,t) \in \mathbb{R} \times [0,T],
\end{equation} with $\CN(u) = |u|^4u$. Let us set
$$I U(x,t) = I\otimes I \otimes I \otimes I(u(x_{1},t)
\otimes u(x_{2},t)\otimes u(x_{3},t)\otimes u(x_{4},t))=
\prod_{j=1}^{4}Iu(x_{j},t).$$ If $u$ solves (\ref{NLSproof}) for
$d=1,$ then $IU$ solves (\ref{NLSproof}) on $\R^4 \times [0,T],$
with right hand side $\CN_I$ given by
$$\CN_I=\sum_{k=1}^{4}(I(\CN_{k})\prod_{j=1,j \ne k}^{4}Iu_{j}).$$
Here and henceforth we set $u_k = u(x_k, t)$, and $\CN_k =
\CN(u_k)$. Hence, according to \eqref{interact3}, we have the
following estimate:

\begin{equation}\label{interactI}
\int_{0}^{T}\int_{\Bbb R^{4} }(-\Delta \Delta a)|IU|^2(x,t)dxdt +2
\int_0^T\int_{\Bbb R^{4}} \nabla a \cdot \{\CN_I,IU\}_pdxdt
\lesssim \sup_{[0,T]}|M_{a}^{I}(t)|,
\end{equation} with $M_{a}^{I}(t)$ the Morawetz action associated
to $IU,$ and $a: \R^4 \rightarrow \R$ convex function.

Now let us decompose,
$$\CN_I=\CN_{good}+\CN_{bad},$$
where
\begin{align*}
\CN_{good} & =\sum_{k=1}^{4}(\CN(Iu_k) \prod_{j=1,j \ne k}^{4}Iu_{j}),\\
\CN_{bad} & = \sum_{k=1}^{4}\left(I(\CN_{k})-\CN(Iu_{k})\right)\prod_{j=1,j \ne k}^{4}Iu_{j}.
\end{align*}
The first summand $\CN_{good}$ creates a defocusing potential and thus
$$\int_{0}^{T}\int_{\Bbb R^4}\nabla a \cdot
\{\CN_{good},IU\}_{p}dxdt \geq 0.$$ Therefore, \eqref{interactI}
yields,

\begin{equation}\label{interactI2}
\int_{0}^{T}\int_{\Bbb R^{4} }(-\Delta \Delta a)|IU|^2(x,t)dxdt
\lesssim \sup_{[0,T]}|M_{a}^{I}(t)| + \left| \int_0^T\int_{\Bbb
R^{4}} \nabla a \cdot \{\CN_{bad},IU\}_pdxdt\right|.
\end{equation}

After performing a change of variables as in Subsection 4.2, and
using the same weight function $a$ we obtain the following
estimate:
\begin{align}\label{amor}
\|Iu\|_{L_{J}^{8}L_{x}^{8}}^{8}
&\lesssim \sup_{[0,T]}\|Iu\|_{\dot{H}^{1}}\|Iu\|_{L^{2}}^{7} \nonumber \\
& + \left|\int_{0}^{T}\int_{\Bbb R^4}\nabla a \cdot
\{\CN_{bad},Iu_1Iu_2Iu_3Iu_4\}_{p}dx_1dx_2dx_3dx_4dt.\right|
\end{align}
Notice that the dot product on the right hand side is left
invariant under the change of variables. Also,
$|\nabla_{z}a(z)|=1$ and hence, since the matrix $A$ is
orthonormal, $|\nabla_{x}a(x)|=1$. Thus we immediately obtain the
$L^\infty$ bound on $\partial_{x_i} a,$ $|\partial_{x_{i}}a(x)|
\leq 1$, $i=1,2,3,4,$ that we strongly use in the following
calculations.

This concludes the proof of the first part of our Theorem. We now
proceed to prove the estimate \eqref{errorcontrol}, which is the
core of this theorem.

\vspace{2mm}

We restrict on a time interval $J=[t_0,t_1]$ on which the solution
$u$ belongs to the space $Z_I(J)$. We wish to compute the dot
product under the sign of integral in \eqref{amor}. First we
observe that $\nabla a$ is real valued, thus
$$\nabla a \cdot \Re (f \overline{\nabla g}-g \overline{\nabla f})=
\Re \left( \nabla a \cdot
(f \overline{\nabla g}-g \overline{\nabla f})\right).$$
Hence, the desired dot product equals
\begin{equation} \label{er1}
\Re\left(\sum_{i=1}^{4}\partial_{x_i} a
\left(\CN_{bad} \; \partial_{x_i}(Iu_1Iu_2Iu_3Iu_4) -
( \prod_{j=1}^{4}Iu_j ) \partial_{x_i}\CN_{bad} \right)\right).
\end{equation}
We start by computing the first summand. Using the definition of
$\CN_{bad},$ and the fact that $\partial_{x_1}$ acts only on $Iu_1$
we obtain the following,
\begin{align*}
&\CN_{bad} \; \partial_{x_1}(Iu_1Iu_2Iu_3Iu_4) -
\left(\prod_{j=1}^{4}Iu_j \right)\partial_{x_1}\CN_{bad}\\
&=\left( \sum_{k=1}^{4}(I(\CN_k) - \CN(Iu_k))\prod_{h=1,h\neq
k}^{4}Iu_h\right)(\partial_{x_1}Iu_1)(Iu_2Iu_3Iu_4)\\
&-\left(\prod_{j=1}^{4}Iu_j \right)\partial_{x_1}\left(
\sum_{k=1}^{4}(I(\CN_k)- \CN(Iu_k))\prod_{h=1,h\neq
k}^{4}Iu_h\right)\\
&=\left( \sum_{k=1}^{4}(I(\CN_k)- \CN(Iu_k))\prod_{h=1,h\neq
k}^{4}Iu_h \right)(\partial_{x_1}Iu_1)(Iu_2Iu_3Iu_4)\\ & -
\prod_{j=1}^{4}Iu_j\left(\sum_{k=1}^{4}\partial_{x_1}(I(\CN_k)-
\CN(Iu_k))\prod_{h=1,h\neq k}^{4}Iu_h\right)\\ &-
\prod_{j=1}^{4}Iu_j\left(\sum_{k=1}^{4}(I(\CN_k)-
\CN(Iu_k))\partial_{x_1}(\prod_{h=1,h\neq k}^{4}Iu_h)\right)\\
&=\left( \sum_{k=1}^{4}(I(\CN_k)- \CN(Iu_k))\prod_{h=1,h\neq
k}^{4}Iu_h \right)(\partial_{x_1}Iu_1)(Iu_2Iu_3Iu_4)\\ & -
\prod_{j=1}^{4}Iu_j \left( \partial_{x_1}(I(\CN_1)-
\CN(Iu_1))\prod_{h=2}^{4}Iu_h \right)\\ &- \prod_{j=1}^{4}Iu_j
\left( \sum_{k=2}^{4}(I(\CN_k)-
\CN(Iu_k))(\partial_{x_1}Iu_1)(\prod_{h=2,h\neq
k}^{4}Iu_h)\right)\\
&=(I(\CN_1)-
\CN(Iu_1))(\partial_{x_1}Iu_1)(Iu_2Iu_3Iu_4)^2\\&+\left( \sum_{k=2}^{4}(I(\CN_k)-
\CN(Iu_k))\prod_{h=1,h\neq
k}^{4}Iu_h\right)(\partial_{x_1}Iu_1)(Iu_2Iu_3Iu_4)\\ & -
\prod_{j=1}^{4}Iu_j\left(\partial_{x_1}(I(\CN_1)-
\CN(Iu_1))\prod_{k=2}^{4}Iu_k\right)\\ &-
\prod_{j=1}^{4}Iu_j\left(\sum_{k=2}^{4}(I(\CN_k)-
\CN(Iu_k))(\partial_{x_1}Iu_1)(\prod_{h=2,h\neq
k}^{4}Iu_h)\right)
\end{align*}
Now notice that the 2nd and 4th
term in the last expression cancel each other, hence
$$\CN_{bad}(\partial_{x_1}(Iu_1Iu_2Iu_3Iu_4)) -
\left(\prod_{j=1}^{4}Iu_j \right)
\partial_{x_1}\CN_{bad}$$
$$=(I(\CN_1)-
\CN(Iu_1))(\partial_{x_1}Iu_1)(Iu_2Iu_3Iu_4)^2$$$$ -
\prod_{j=1}^{4}Iu_j\left(\partial_{x_1}(I(\CN_1)-
\CN(Iu_1))\prod_{k=2}^{4}Iu_k\right)$$

$$=\left[(I(\CN_1)-
\CN(Iu_1))\partial_{x_1}Iu_1-\partial_{x_1}(I(\CN_1)-
\CN(Iu_1))Iu_1\right](Iu_2Iu_3Iu_4)^2.
$$

\

\noindent Hence the first summand in \eqref{er1} is given
by,\vspace{1mm}
$$\Re \left( \partial_{x_1}a\left[(I(\CN_1)-
\CN(Iu_1))\partial_{x_1}Iu_1-\partial_{x_1}(I(\CN_1)-
\CN(Iu_1))Iu_1\right](Iu_2Iu_3Iu_4)^2 \right).$$

\noindent Analogously one can see that the $i$th summand,
$i=1,2,3,4$ is given by:
$$\Re \left( \partial_{x_i}a\left[(I(\CN_i)-
\CN(Iu_i))\partial_{x_i}Iu_i-\partial_{x_i}(I(\CN_i)-
\CN(Iu_i))Iu_i\right](\prod_{j=1,j\neq i}^{j=4}Iu_j)^2 \right).$$
Thus, our error term
\begin{align}\mathcal{E}=\int_{t_0}^{t_1}\int_{\mathbb{R}^4}\nabla & a \cdot \left(\CN_{bad}(\nabla(Iu_1Iu_2Iu_3Iu_4)) -
\prod_{j=1}^{4}Iu_j\nabla \CN_{bad}\right)\nonumber\end{align}
reduces to
\begin{align*}
\mathcal{E}=\Re ( \int_{t_0}^{t_1}\int_{\mathbb{R}^4}\sum_{k=1}^{4}&\{\partial_{x_k}a\left[(I(\CN_k)-
\CN(Iu_k))\partial_{x_k}Iu_k-\partial_{x_k}(I(\CN_k)-
\CN(Iu_k))Iu_k\right]\\&\times(\prod_{j=1,j\neq
k}^{4}Iu_j)^2\}dx_1dx_2dx_3dx_4dt ).
\end{align*}

\noindent Hence, by symmetry,
\begin{equation}\label{total error}|\mathcal{E}| \lesssim |E|,
\end{equation}
where
\begin{align*}E=\int_{t_0}^{t_1}\int_{\mathbb{R}^4}&\{\partial_{x_1}a\left[(I(\CN_1)-
\CN(Iu_1))\partial_{x_1}Iu_1-\partial_{x_1}(I(\CN_1)-
\CN(Iu_1))Iu_1\right]\\&\times(\prod_{j=2}^{4}Iu_j)^2\}dx_1dx_2dx_3dx_4dt.\end{align*}

We have,

\begin{equation}\label{E}|E| \leq E_1 + E_2\end{equation} where

$$E_1= \int_{t_0}^{t_1}\int_{\mathbb{R}^4}|\partial_{x_1}a||I(\CN_1)-
\CN(Iu_1)||\partial_{x_1}Iu_1|\prod_{j=2}^{4}|Iu_j|^2dx_1dx_2dx_3dx_4dt$$
and

$$E_2= \int_{t_0}^{t_1}\int_{\mathbb{R}^4}|\partial_{x_1}a||\partial_{x_1}(I(\CN_1)-
\CN(Iu_1))||Iu_1|\prod_{j=2}^{4}|Iu_j|^2dx_1dx_2dx_3dx_4dt.$$
Since $|\partial_{x_1}a| \leq 1$, after applying Fubini's
theorem we have

$$E_1 \leq \left(\int_{t_0}^{t_1}\int_{\mathbb{R}}|I(\CN_1)-
\CN(Iu_1)||\partial_{x_1}Iu_1|dx_1dt\right)\|Iu\|^6_{L^{\infty}_tL^2_{x}},$$
and

$$E_2 \leq \left(\int_{t_0}^{t_1}\int_{\mathbb{R}}|\partial_{x_1}(I(\CN_1)-
\CN(Iu_1))||Iu_1|dx_1dt\right)\|Iu\|^6_{L^{\infty}_tL^2_{x}}.$$
Since the pair $(\infty, 2)$ is admissible, we then obtain (rename
$x_1=x$ ):

$$E_1 \leq \left(\int_{t_0}^{t_1}\int_{\mathbb{R}}|I(\CN)-
\CN(Iu)||\partial_{x}Iu|dx\; dt\right) \|u\|^6_{Z_I(J)},$$ and

$$E_2 \leq \left(\int_{t_0}^{t_1}\int_{\mathbb{R}}|\partial_{x}(I(\CN)-
\CN(Iu))||Iu|dx\; dt\right)\|u\|^6_{Z_I(J)}.$$ Therefore,

$$E_1 \leq \|I(\CN)-
\CN(Iu)\|_{L^1_tL^2_x}\|\partial_{x}Iu\|_{L^\infty_tL^2_x}\;
\|u\|^6_{Z_I(J)},$$ and

$$E_2 \leq \|\partial_x(I(\CN)-
\CN(Iu))\|_{L^1_tL^2_x}\|Iu\|_{L^\infty_tL^2_x}\;\|u\|^6_{Z_I(J)}.$$
Again, since $(\infty,2)$ is admissible we obtain:

$$E_1 \leq \|I(\CN)-
\CN(Iu)\|_{L^1_tL^2_x}\; \|u\|^7_{Z_I(J)},$$ and

$$E_2 \leq \|\partial_x(I(\CN)-
\CN(Iu))\|_{L^1_tL^2_x}\; \|u\|^7_{Z_I(J)}.$$ Therefore, from
\eqref{E} and the bounds above, we deduce that
\begin{equation}\label{error}|E| \leq \left( \|I(\CN)- \CN(Iu)\|_{L^1_tL^2_x}+\|\partial_x(I(\CN)-
\CN(Iu))\|_{L^1_tL^2_x}\right)\|u\|^7_{Z_I(J)}.\end{equation}

We proceed to estimate $\|\partial_x(I(\CN)- \CN(Iu))\|_{L^1_tL^2_x}$,
which is the hardest of the two terms. Toward this aim, let us
observe that since $\CN(u)=|u|^pu$ with $p=4$, we will be able to
work on the Fourier side to estimate the commutator $I(\CN)- \CN(Iu)$.

We compute,\footnote{We ignore complex conjugates, since our
computations are not effected by conjugation.}

$$\widehat{\partial_x(I(\CN)- \CN(Iu))}(\xi) = \int_{\xi=\xi_1+...+\xi_5} i\xi[m(\xi)-m(\xi_1)\cdots m(\xi_5)]
\hat{u}(\xi_1) \ldots \hat{u}(\xi_5) d\xi_1 \ldots d\xi_5.$$ We
decompose $u$ into a sum of dyadic pieces  localized around $N_j$
in the usual way. In all the estimates that follow we obtain a factor of order $N_{max}^{-\epsilon}$ in order
 to be able to perform the summations at the end. We omit this technical detail. Then,

\begin{align} & \|\partial_x(I(\CN)-
\CN(Iu))\|_{L^1_tL^2_x}=\|\widehat{\partial_x(I(\CN)-
\CN(Iu))}\|_{L^1_tL^2_\xi}\label{fourierside}\\ \nonumber&\leq
\sum_{N_1,\ldots, N_5} \|\int_{\xi=\xi_1+...+\xi_5; |\xi_i|\sim
N_i} \xi[m(\xi)-m(\xi_1)\cdots m(\xi_5)] \widehat{u_1}\ldots
\widehat{u_5} d\xi_1 \ldots d\xi_5\|_{L^1_tL^2_\xi}\\
\nonumber&=\sum_{N_1,\ldots, N_5} \|\int_{\xi=\xi_1+...+\xi_5;
|\xi_i|\sim N_i} \xi\frac{[m(\xi)-m(\xi_1)\cdots
m(\xi_5)]}{m(\xi_1)\cdots m(\xi_5)} \widehat{Iu_1}\ldots
\widehat{Iu_5} d\xi_1 \ldots
d\xi_5\|_{L^1_tL^2_\xi}.\end{align}

Without loss of generality, we can assume that the $N_j$'s are
rearranged so that $$N_1 \geq \ldots \geq N_5.$$ Set,

$$\sigma(\xi_1,\ldots, \xi_5)=(\xi_1 + \ldots +\xi_5)\frac{[m(\xi_1 + \ldots +\xi_5)-m(\xi_1)\cdots m(\xi_5)]}
{m(\xi_1)\cdots m(\xi_5)}.$$
Then,
$$\sigma(\xi_1,\ldots,\xi_5)= \sum_{j=1}^{6}\sigma_j(\xi_1,\ldots,\xi_5),$$
with
$$\sigma_j(\xi_1,\ldots,\xi_5) = \chi_j(\xi_1,\ldots,\xi_5) \sigma(\xi_1,\ldots,\xi_5),$$
where
$\chi_j(\xi_1,\ldots,\xi_5)$ are  smooth characteristic functions
of the sets $\Omega_j$ defined as follows:
\begin{itemize}
\item $\Omega_1=\{|\xi_i| \sim N_i, i=1,\ldots,5 ; N_1 \ll N\}$.
\item$\Omega_2=\{|\xi_i| \sim N_i, i=1,\ldots,5 ;N_1 \gtrsim N \gg
N_2\}.$ \item $\Omega_3=\{|\xi_i| \sim N_i, i=1,\ldots,5 ;N_1 \geq
N_2 \gtrsim N \gg N_3.\}$ \item $\Omega_4=\{|\xi_i| \sim N_i,
i=1,\ldots,5 ;N_1 \geq N_2 \geq N_3 \gtrsim N \gg N_4\}.$ \item
$\Omega_5=\{|\xi_i| \sim N_i, i=1,\ldots,5 ;N_1 \geq N_2 \geq N_3
\geq N_4 \gtrsim N \gg N_5\}.$ \item $\Omega_6=\{|\xi_i| \sim N_i,
i=1,\ldots,5 ;N_1, \ldots, N_5 \gtrsim N\}.$
\end{itemize}
Hence, from \eqref{fourierside} we get,
\begin{align}&\label{bound}\|\partial_x(I(\CN)-
\CN(Iu))\|_{L^1_tL^2_x}\\ &\lesssim \sum_{N_1,\ldots, N_5}
\sum_{j=1}^{6}\|\int_{\xi=\xi_1+...+\xi_5}
\sigma_j(\xi_1,\ldots,\xi_5)\widehat{Iu_1}\ldots
\widehat{Iu_5} d\xi_1 \ldots d\xi_5\|_{L^1_tL^2_\xi}
=\sum_{N_1,\ldots, N_5} \sum_{j=1}^{6}L_j\nonumber.\end{align}

We proceed to analyze the contribution of each of the integrals
$L_j.$

\

\textbf{Contribution of $L_1$.} Since $\sigma_1$ is identically
zero, $L_1$ gives no contribution to the sum above.

\

\textbf{Contribution of $L_2$.} We have,

$$\|\int_{\xi=\xi_1+...+\xi_5}
\sigma_2(\xi_1+\ldots+\xi_5) \widehat{Iu_1}\ldots
\widehat{Iu_5} d\xi_1 \ldots d\xi_5\|_{L^1_tL^2_\xi}
$$
$$=\frac{1}{N}\|\int_{\xi=\xi_1+...+\xi_5}
\frac{N}{\xi_1\xi_2}\sigma_2(\xi_1,\ldots,\xi_5)\widehat{\langle \partial_{x} \rangle Iu_1}
\widehat{\langle \partial_{x} \rangle Iu_2}\ldots \widehat{Iu_5} \, d\xi_1
\ldots d\xi_5\|_{L^1_tL^2_\xi}$$
$$\lesssim \frac{1}{N} \|\langle \partial_{x} \rangle Iu_1\|_{L^5_tL^{10}_x}
\|\langle\partial_{x}\rangle Iu_2\|_{L^5_tL^{10}_x}\prod_{j=3}^{5}\|Iu_j\|_{L^5_tL^{10}_x}$$
where in the last line we used the Coifman-Meyer multiplier
theorem, and H\"{o}lder in time. The application of the multiplier
theorem is justified by the fact that the symbol
$$a_2(\xi_1,\ldots,\xi_5)=
\frac{N}{\xi_1\xi_2}\sigma_2(\xi_1,\ldots,\xi_5)$$ is of order
zero. The $L^\infty$ bound follows after an application of the
mean value theorem. Indeed,

$$|a_2(\xi_1,...,\xi_5)| \leq \frac{N}{N_1N_2}|\xi_1+\ldots+\xi_5|\frac{|\nabla_{\xi_{1}} m(\xi_1)(\xi_2+\ldots+\xi_5)
|}{m(\xi_1)} \lesssim N_1 \frac{N}{N_1N_2} \frac{N_2}{N_1}\lesssim
1.$$

\

\textbf{Contribution of $L_3$.} We have,

$$\|\int_{\xi=\xi_1+...+\xi_5}
\sigma_3(\xi_1+\ldots+\xi_5) \widehat{Iu_1}\ldots
\widehat{Iu_5} \, d\xi_1 \ldots d\xi_5\|_{L^1_tL^2_\xi}
$$
$$=\frac{1}{N}\|\int_{\xi=\xi_1+...+\xi_5}
\frac{N}{\xi_1\xi_2}\sigma_3(\xi_1,\ldots,\xi_5)\widehat{\langle \partial_{x} \rangle
Iu_1}\widehat{\langle \partial_{x} \rangle Iu_2}\ldots \widehat{Iu_5} \, d\xi_1
\ldots d\xi_5\|_{L^1_tL^2_\xi}$$
$$\lesssim \frac{1}{N} \|\langle\partial_{x}\rangle Iu_1\|_{L^5_tL^{10}_x}
\|\langle \partial_{x} \rangle Iu_2\|_{L^5_tL^{10}_x}\prod_{j=3}^{5}\|Iu_j\|_{L^5_tL^{10}_x}$$
where in the last line we used the Coifman-Meyer multiplier
theorem, and Holder in time. The application of the multiplier
theorem is justified by the fact that the symbol
$$a_3(\xi_1,\ldots,\xi_5)=
\frac{N}{\xi_1\xi_2}\sigma_3(\xi_1,\ldots,\xi_5)$$ is of order
zero. The $L^\infty$ bound follows from the following chain of
inequalities,

$$|a_3(\xi_1, \ldots, \xi_5)| \lesssim \frac{N}{N_1N_2}
|\xi_1+\ldots+\xi_5||(\frac{m(\xi_1+\ldots+\xi_5)}{m(\xi_1)m(\xi_2)}+1)
$$$$ \lesssim \frac{N}{N_1N_2}(\frac{N_1}{m(N_2)} + N_1) \lesssim 1.$$
We have used the fact that $|\xi|m(\xi)$ is monotone increasing and thus
$$|(\xi_1+\ldots+\xi_5)m(\xi_1+\ldots+\xi_5)| \lesssim N_{1}m(\xi_{1}).$$

It is now evident what is the contribution of the remaining cases.

\

\textbf{Contribution of $L_4$.}

$$\|\int_{\xi=\xi_1+...+\xi_5}
\sigma_4(\xi_1+\ldots+\xi_5) \widehat{Iu_1}\ldots
\widehat{Iu_5} \, d\xi_1 \ldots d\xi_5\|_{L^1_tL^2_\xi}
$$$$\lesssim \frac{1}{N^2} \|\langle \partial_{x}\rangle
Iu_1\|_{L^5_tL^{10}_x}\|\langle \partial_{x} \rangle
Iu_2\|_{L^5_tL^{10}_x}\|\langle \partial_{x} \rangle
Iu_3\|_{L^5_tL^{10}_x}\prod_{j=4}^{5}\|Iu_j\|_{L^5_tL^{10}_x}$$

where in this case the symbol to which we apply the multiplier
theorem is:

$$a_4(\xi_1,\ldots,\xi_5)=
\frac{N^2}{\xi_1\xi_2\xi_3}\sigma_4(\xi_1,\ldots,\xi_5).$$

\

\textbf{Contribution of $L_5$.}
$$\|\int_{\xi=\xi_1+...+\xi_5}
\sigma_5(\xi_1+\ldots+\xi_5) \widehat{Iu_1}\ldots
\widehat{Iu_5} d\xi_1 \ldots d\xi_5\|_{L^1_tL^2_\xi}
$$$$\lesssim \frac{1}{N^3} \prod_{j=1}^{4}\|\langle \partial_{x} \rangle
Iu_j\|_{L^5_tL^{10}_x}\|Iu_5\|_{L^5_tL^{10}_x}$$ where in this
case the symbol to which we apply the multiplier theorem is:

$$a_5(\xi_1,\ldots,\xi_5)=
\frac{N^3}{\xi_1\xi_2\xi_3\xi_4}\sigma(\xi_1,\ldots,\xi_5).$$

\

\textbf{Contribution of $L_6$.}

$$\|\int_{\xi=\xi_1+...+\xi_5}
\sigma_6(\xi_1+\ldots+\xi_5) \widehat{Iu_1}\ldots
\widehat{Iu_5} d\xi_1 \ldots d\xi_5\|_{L^1_tL^2_\xi}\lesssim
\frac{1}{N^4} \prod_{j=1}^{5}\|\langle \partial_{x} \rangle Iu_j\|_{L^5_tL^{10}_x},$$
where in this case the symbol to which we apply the multiplier
theorem is:

$$a_6(\xi_1,\ldots,\xi_5)=
\frac{N^4}{\xi_1\xi_2\xi_3\xi_4\xi_5}\sigma_6(\xi_1,\ldots,\xi_5).$$

\

In all the cases above, we proved the $L^\infty$ bound for the
symbols $a_i(\xi_1,\ldots,\xi_5), i=2,\ldots,6.$ We proceed here
to determine the bound on $\partial_{\xi_1} a_2.$ The remaining
bounds are left to the reader. Recall
that$$a_2(\xi_1,\ldots,\xi_5)=
\frac{N}{\xi_1\xi_2}\sigma(\xi_1,\ldots,\xi_5)\chi_2(\xi_1,\ldots,\xi_5).$$
Hence,

$$|\partial_{\xi_1}a_2| \lesssim |\partial_{\xi_1}(\frac{N}{\xi_1\xi_2}\sigma(\xi_1,\ldots,\xi_5))
\chi_2(\xi_1,\ldots,\xi_5)|+$$
$$|\partial_{\xi_1}(\chi_2(\xi_1,\ldots,\xi_5))\frac{N}{\xi_1\xi_2}\sigma(\xi_1,\ldots,\xi_5)|.$$

The bound on the second summand, follows as the $L^\infty$ bound
on $a_2.$ We proceed to bound the first summand, that in turn is
bounded by the sum of the following two terms:

\begin{equation}\label{firstterm}|N \frac{\xi_1\xi_2 -
\xi_2(\xi_1+\ldots+\xi_5)}{\xi_1^2\xi_2^2}[\frac{m(\xi_1+\ldots+\xi_5)
- m(\xi_1)}{m(\xi_1)}]|\end{equation}

\begin{equation}\label{secondterm}|N \frac{
\xi_1+\ldots+\xi_5}{\xi_1\xi_2}[\frac{\partial_{\xi_1}m(\xi_1+\ldots+\xi_5)}{m(\xi_1)}
- m(\xi_1+\ldots+\xi_5)\frac{\partial_{\xi_1}m(\xi_1)}{m^2(\xi_1)}
]|.\end{equation} Again, an application of the mean value theorem
gives that

$$\eqref{firstterm} \lesssim \frac{N}{N_1N_2}\frac{N_2}{N_1} \lesssim 1.$$
As for \eqref{secondterm}, it is easy to see that

$$|(\xi_1+\ldots+\xi_5)\frac{\partial_{\xi_1}m(\xi_1+\ldots+\xi_5)}{m(\xi_1)}| \lesssim 1$$
while, using also the monotonicity of $|\xi|m(\xi)$ in the form
$$|(\xi_1+\ldots+\xi_5)m(\xi_1+\ldots+\xi_5)| \lesssim N_{1}m(\xi_{1})$$
one gets
$$|(\xi_1+\ldots+\xi_5)m(\xi_1+\ldots+\xi_5)\frac{\partial_{\xi_1}m(\xi_1)}{m^2(\xi_1)}| \lesssim 1.$$
Thus,
$$\eqref{secondterm} \lesssim 1,$$ and we obtain the desired bound
on $\partial_{\xi_1}a_2.$
\\
\\
Finally, since the pair $(5,10)$ is admissible, we obtain that in
all of the cases above

$$\|\int_{\xi=\xi_1+...+\xi_5}
\sigma_i(\xi_1+\ldots+\xi_5) \widehat{Iu_1}\ldots
\widehat{Iu_5} d\xi_1 \ldots d\xi_5\|_{L^1_tL^2_\xi}
$$$$\lesssim \frac{1}{N}\|u\|^5_{Z_I(J)}.$$
Therefore, we deduce from \eqref{bound} that
$$\|\partial_x(I(\CN)- \CN(Iu))\|_{L^1_tL^2_x} \lesssim
\frac{1}{N^-}\|u\|^5_{Z_I(J)}.$$ Analogously, $$\|I(\CN)-
\CN(Iu)\|_{L^1_tL^2_x} \lesssim \frac{1}{N^-}\|u\|^5_{Z_I(J)}.$$
Hence, in view of \eqref{error} we obtain the following estimate
for the error term,

$$|E| \lesssim \frac{1}{N^-}\|u\|^{12}_{Z_I(J)}.$$
Thus, \eqref{total error} implies
$$\left|\int_{t_0}^{t_1}\int_{\mathbb{R}^4}\nabla a \cdot \left(\CN_{bad}(\nabla(Iu_1Iu_2Iu_3Iu_4)) -
\prod_{j=1}^{4}Iu_j\nabla \CN_{bad}\right)\right| \lesssim
\frac{1}{N^-}\|u\|^{12}_{Z_I(J)},$$ which concludes the proof.
\\

\section{Proof of Theorem \ref{th-gwp}} \label{sec-proof}

Suppose that $u(x,t)$ is a global in time solution to \eqref{ivp1}
with  initial data $u_0 \in C_0^\infty(\R^n)$. Set
$u^\lambda(x)=\frac{1}{\lambda^{\frac{1}{2}}}u(\frac{x}{\lambda},\frac{t}{\lambda^{2}})$.
We choose the parameter $\lambda$ so that $\|I u_0^\lambda\|_{H^1}
= O(1)$, that is
\begin{equation*}\label{L}
\lambda \sim N^{\frac{1-s}{s}}.
\end{equation*}
Next, let us pick a time
$T_0$ arbitrarily large, and let us define
\begin{equation*}
S : = \{0 < t < \lambda^2T_0 :
\|Iu^\lambda\|_{L^{6}_tL^{6}_x([0,t]\times
\R )} \leq Kt^{\frac{1}{18}}N^{\frac{1}{9}}\},
\end{equation*}
with $K$ a constant to be chosen later. We claim that $S$ is the
whole interval $[0,\lambda^2T_0].$ Indeed, assume by contradiction
that it is not so, then since
$$\|Iu^\lambda\|_{L^{6}_tL^{6}_x([0,t]\times
\R )}$$ is a continuous function of time, there exists a time $T
\in [0,\lambda^2T_0]$ such that
\begin{align}
\label{contr1}
& \|Iu^\lambda\|_{L^{6}_tL^{6}_x([0,T]\times
\R )} >KT^{\frac{1}{18}}N^{\frac{1}{9}}\\
\label{contr2}
&\|Iu^\lambda\|_{L^{6}_tL^{6}_x([0,T]\times
\R )} \leq 2KT^{\frac{1}{18}}N^{\frac{1}{9}}.
\end{align}
We now split the interval $[0,T]$ into subintervals
$J_k$, $k=1,...,L$ in such a way
that
\begin{equation*}
\|Iu^\lambda\|_{L^{6}_tL^{6}_x(J_{k}\times\R )}^6 \leq \mu,
\end{equation*}
with $\mu$ as in Proposition \ref{lwp}. This is possible because of \eqref{contr2}. Then, the
number $L$ of possible subintervals must satisfy
\begin{equation} \label{L}
L \sim \frac{(2KT^{\frac{1}{18}}N^{\frac{1}{9}})^{6}}{\mu} \sim \frac{(2K)^{6}T^{\frac{1}{3}}N^{\frac{2}{3}}}{\mu}.
\end{equation}
From Proposition \ref{lwp} and Propositions \ref{e1} and \ref{e2} we know that,
for any $1/3<s<1$
\begin{equation*}\label{energybound}
\sup_{[0,T]}E(Iu^\lambda(t)) \lesssim E(Iu_0^\lambda) +
\frac{L}{N^{2}}
\end{equation*}
and by our choice of $\lambda$, $E(Iu_0^\lambda)\lesssim 1.$
Note that if we restrict to $s>1/3$ we can apply the previous Propositions.
Hence, in order to guarantee that
\begin{equation*}\label{energyb}
E(Iu^\lambda)\lesssim 1\end{equation*} holds for all $t \in [0,T]$
we need to require  that
$$ L \lesssim N^{2}. $$
Since $T \leq \lambda^2T_0,$ according to \eqref{L}, this is fulfilled
as long as
\begin{equation} \label{LN}
\frac{(2K)^{6} (\lambda^2 T_0)^{\frac{1}{3}}N^{\frac{2}{3}}}{\mu} \sim
N^{2}.
\end{equation}
From our choice of $\lambda$, the expression \eqref{LN} implies that
\begin{equation*}\label{Tn}
T_{0}^{\frac{1}{3}}\frac{(2K)^6}{\mu} \sim N^{\frac{4}{3}-\frac{2(1-s)}{3s}}
= N^{\frac{2(3s-1)}{3s}}.
\end{equation*}
Thus if $s> 1/3$, we have that $N$ is a large number for $T_{0}$ large.
\\
\\
Now recall the a priori estimate \eqref{amor}
\begin{align*}
\|Iu\|_{L_{T}^{8}L_{x}^{8}}^{8} &\lesssim
\sup_{0,T}\|Iu\|_{\dot{H}^{1}}\|Iu\|_{L^{2}}^{7}+\\\nonumber
&\left|\int_0^T\int_{\Bbb R^4}\nabla a \cdot
\{N_{bad},Iu_1Iu_2Iu_3Iu_4\}_{p}dx_1dx_2dx_3dx_4dt\right|.
\end{align*}
Set
$$Error(t):=\int_{\Bbb R^4}\nabla a \cdot
\{N_{bad},Iu_1Iu_2Iu_3Iu_4\}_{p}dx_1dx_2dx_3dx_4.$$ By Theorem
\ref{main} and Proposition \ref{lwp} on each interval $J_{k}$ we
have that
$$\left|\int_{J_{k}}Error(t)dt\right| \lesssim \frac{1}{N}\|u\|^{12}_{Z_{I}}
 \lesssim \frac{1}{N}\|I
u^\lambda(t_{0})\|_{H^1}^{12} \lesssim \frac{1}{N}.$$
Summing all the $J_{k}$'s we have that
$$\left|\int_{0}^{T}Error(t)dt\right|
\lesssim L\frac{1}{N} \sim \frac{N^{2}}{N} \sim
N.$$
Therefore,
$$\|Iu^{\lambda}\|_{L_{T}^{8}L_{x}^{8}}^{8} \lesssim
\sup_{t \in [0,T]}\|Iu^{\lambda}\|_{\dot{H}^{1}}\|Iu^{\lambda}\|_{L^{2}}^{7}
+\left|\int_{0}^{T}Error(t)dt\right|
\lesssim 1+N \sim N,$$ which implies
\begin{equation}\label{l8}
\|Iu^{\lambda}\|_{L_{T}^{8}L_{x}^{8}}\lesssim N^{\frac{1}{8}}.
\end{equation}
On the other hand H\"older inequality in time together with the definition of the $I$ operator and conservation of mass gives
\begin{equation}\label{l2}
\|Iu^{\lambda}\|_{L_{t \in [0,T]}^{2}L_{x}^{2}}
\leq T^{\frac{1}{2}}\|Iu^{\lambda}\|_{L_{t \in [0,T]}^{\infty}L_{x}^{2}}
\lesssim T^{\frac{1}{2}}\|u^{\lambda}\|_{L_{x}^2}=T^{\frac{1}{2}}\|u_{0}\|_{L_{x}^2}\sim
T^{\frac{1}{2}}.
\end{equation}
Interpolation between \eqref{l8} and \eqref{l2} gives that
$$\|Iu^{\lambda}\|_{L_{t \in [0,T]}^{6}L_{x}^{6}} \leq CT^{\frac{1}{18}}N^{\frac{1}{9}}.$$
This estimate contradicts \eqref{contr1} for an appropriate
choice of $K$. Hence $S = [0,\lambda^2T_0]$, and $T_0$ can
be chosen arbitrarily large. In addition, we have also proved that for $s>1/3$
$$\| I u^\lambda(\lambda^{2}T_{0})\|_{H^1_x}=O(1).$$
Then,
$$\|u(T_{0})\|_{H^{s}} \lesssim \|u(T_{0})\|_{L^{2}}+\|u(T_{0})\|_{\dot{H}^{s}}=
\|u_{0}\|_{L^{2}}+\lambda^{s}\|u^{\lambda}(\lambda^{2}T_{0})\|_{\dot{H}^{s}}$$
$$\lesssim \lambda^{s}\| I u^\lambda(\lambda^{2}T_{0})\|_{H^1_x}\lesssim \lambda^{s}\lesssim N^{1-s}
\lesssim T_{0}^{\frac{s(1-s)}{2(3s-1)}}.$$ Since $T_{0}$ is
arbitrarily
 large, the a priori bound on the $H^{s}$ norm concludes the global well-posedness of the the Cauchy
problem \eqref{ivp1}-\eqref{bc1}.

\end{document}